\documentclass[11pt]{article}
\usepackage[a4paper, total={6in, 8in}]{geometry}
\usepackage{authblk}
\usepackage{amsmath}
\usepackage{amsthm}
\usepackage[OT1]{fontenc}
\numberwithin{equation}{section}
\usepackage{amsfonts}
\usepackage{amssymb}
\usepackage{graphicx}
\usepackage{esint}
\usepackage{hyperref}
\usepackage{bm}
\usepackage{comment}
\usepackage{tikz-cd}
\hypersetup{linktocpage,
colorlinks=true,
linktoc=all,
linkcolor=blue,
}

\newcommand{\nd}[2]{\makebox[7em]{\begin{tabular}{@{}c@{}}#1\\[-1pt]\scriptsize #2\end{tabular}}}
\newcommand{\aname}[1]{\text{\footnotesize #1}}
\newcommand{\lref}[1]{\text{\scriptsize #1}}
\newcommand{\twolab}[2]{\begin{tabular}{@{}c@{}}\footnotesize #1\\[-1pt]\scriptsize #2\end{tabular}}

\makeatletter

\newtheorem{theorem}{Theorem}[section]

\newtheorem{definition}{Definition}[section]

\newtheorem{lemma}{Lemma}[section]
\newtheorem{proposition}{Proposition}[section]
\newtheorem{remark}{Remark}[section]

\usepackage[color=magenta]{todonotes}

 \newcommand{\D}{{\mathcal{D}}}
 \newcommand{\HS}{{\mathcal{D}^{\rm hs}}}

\newcommand{\DsL}{{\mathcal{D}}^{{\rm st}, L}} 
\newcommand{\DhL}{{\mathcal{D}}^{{\rm hs}, L}} 
\newcommand{\Dis}{{\mathcal{D}}^{\rm st}} 
\newcommand{\Dh}{{\mathcal{D}}^{\rm hs}} 

\newcommand{\R}{{\mathbb{R}}}

\newcommand{\Ra}{{\rm{Ra}}}
\newcommand{\Nu}{{\rm{Nu}}}

\newcommand{\ub}{{\bm{u}}}
\newcommand{\xx}{{\bm{x}}}
\newcommand{\yy}{{\bm{y}}}
\newcommand{\zz}{{\bm{z}}}
\newcommand{\kk}{{\bm{k}}}
\newcommand{\el}{\bm{\ell}}
\newcommand{\vv}{{\bm{v}}}

\newcommand{\xib}{{\bm{\xi}}}
\newcommand{\etab}{{\bm{\eta}}}

\newcommand{\QQ}{\mathcal{Q}}
\newcommand{\ZZ}{\mathbb{Z}}
\newcommand{\TT}{\mathbb{T}}

\newcommand{\bzero}{{\bm{0}}}

\newcommand{\rr}{{\bm{r}}}

\title{Covariance Structures and Transport Regimes in \\ Stochastic Convection}

\author[1]{Franco Flandoli}
\author[2]{Theresa Lange}
\author[3]{Camilla Nobili}
\affil[1]{Scuola Normale Superiore di Pisa, Piazza dei Cavalieri 7, 56126 Pisa, Italy}
\affil[2]{Max-Planck Institute for Mathematics in the Sciences, Inselstr. 22-26, 04103 Leipzig, Germany}
\affil[3]{School of Mathematics and Physics, University of Surrey, Guildford, Surrey GU2 7HX, UK}
\affil[ ]{E-mail: franco.flandoli@sns.it, theresa.lange@mis.mpg.de, c.nobili@surrey.ac.uk}
\begin{document}

\maketitle
\begin{abstract}
We study the scaling of the Nusselt number $\Nu$ in wall-bounded thermal convection at infinite Prandtl number from a new, physically motivated point of view. Within a stochastic framework, instead of forcing the system externally, we model the internal temperature fluctuations by a noise that is white in time and coloured in space: these fluctuations drive the velocity field through the Stokes equations, and the resulting velocity field transports the temperature. The noise is therefore not an external perturbation, but a parametrisation of the turbulent fluctuations of the system itself, in the spirit of decoupling approaches to turbulent transport. Within this model, we identify conditions on the spatial covariance $q_S$ of the fluctuations under which the two prominent scaling laws are recovered: if the horizontal correlation length is of order one up to the plates, we obtain the Kraichnan--Spiegel bound $\Nu\lesssim\Ra^{1/2}$; if, instead, the fluctuations have a boundary layer of thickness $\eta$, possibly depending on $\Ra$, in which the correlation length shrinks to zero at the plates, we obtain
$\Nu\lesssim\Ra^{1/3}+(\Ra\,\eta)^{1/2}$; in particular, the Malkus scaling $\Ra^{1/3}$ holds as soon as this boundary layer is no thicker than the thermal one, $\eta\lesssim\Ra^{-1/3}$. All bounds are uniform in the aspect ratio. The proof combines sharp estimates on the explicit half-space kernel with an energy estimate for the correction induced by the upper plate.
\end{abstract}
\newpage 
\tableofcontents
\section{Introduction}
The aim of this paper is to understand how the statistics of temperature fluctuations near the
boundaries determine the scaling of the Nusselt number $\Nu$ in strongly convective regimes. We
consider a model for a fluid confined between two horizontally aligned plates in which the velocity field is driven by a stochastic temperature field rather
than by the temperature itself, and we study under which conditions on the statistics of this field
the upper bounds on heat transport known for Rayleigh-B\'enard convection at infinite Prandtl
number are recovered. The model deviates from the Rayleigh-B\'enard system, but it retains the
same transport mechanism.

\smallskip

In thermal convection, the velocity field $\ub$ is driven by buoyancy forces, which are caused by
temperature variations along the direction of gravity. When the difference between the temperatures imposed at the
two plates is large, the fluid becomes turbulent. A prototypical example is the three-dimensional
Rayleigh--B\'enard model at infinite Prandtl number. There, the temperature $T$ is advected by $\ub$
and diffuses,
\[
\partial_tT+\ub\cdot\nabla T=\Delta T,
\]
with $T=1$ at the bottom plate and $T=0$ at the top plate, and the velocity is coupled to the
temperature through the Stokes equations
\begin{equation}\label{eq:UStokes}
-\Delta\ub+\nabla p=\Ra\,T\,\bm{e}_{x_3},\qquad\nabla\cdot\ub=0,
\end{equation}
with no-slip boundary conditions at both plates; here $p$ denotes the pressure and $\bm{e}_{x_3}$ the unit vector in the vertical direction $x_3$. In the original
dimensional formulation, the temperatures at the bottom and top plates are $T_{\rm bottom}$
and $T_{\rm top}$, and after non-dimensionalisation the temperature difference
$T_{\rm bottom}-T_{\rm top}$ enters only through the Rayleigh number
\[
\Ra=\frac{\alpha g(T_{\rm bottom}-T_{\rm top})h^3}{\nu\kappa}.
\]
Here $g$ is the gravitational acceleration, $\alpha$ the thermal expansion coefficient, $h$ the
height of the layer, and $\nu$ and $\kappa$ the kinematic viscosity and the thermal diffusivity, respectively.
The turbulent regime corresponds to the limit $\Ra\to\infty$.

The main quantity of interest is the Nusselt number $\Nu$, the averaged vertical heat flux
normalised by the purely conductive one. A central question \cite{AGL09} is how $\Nu$ scales with
$\Ra$, with the Prandtl number $\Pr=\nu/\kappa$, and with the aspect ratio $\Gamma$ of the
container (denoted by $L$ below); from the mathematical point of view, one would like to prove
rigorous upper bounds of the form
\[
\Nu\le F(\Ra,\Pr,\Gamma)
\]
with an explicit function $F$. Within the ongoing debate about the \textit{ultimate scaling} of turbulent thermal convection \cite{lohse2024ultimate,doering2019}, two scaling laws play a central role: the
\emph{Kraichnan--Spiegel scaling} \cite{K1962,Sp1963}
\[
\Nu\sim\Ra^{1/2},
\]
and the \emph{Malkus scaling} \cite{Malkus1954}
\[
\Nu\sim\Ra^{1/3}.
\]
The bound $\Nu\lesssim\Ra^{1/2}$ holds rigorously for every Prandtl number \cite{howard1963heat,DC96}.
At infinite Prandtl number, the Malkus scaling was first proved up to logarithmic corrections
\cite{CD99,DOR06,otto2011rayleigh}. Most recently Chanillo and Malchiodi 
\cite{ChanMal2025}, building up on the direct PDE method of Otto \& al. \cite{otto2011rayleigh, choffrutNobiliOttoUpperBoundsOnNusseltNumberAtFinitePrandtlNumber, NO17} succeeded in removing logarithmic corrections performing a sharp analysis on the integral representation of the fourth-order boundary value problem for the vertical velocity $u_3$, using a crucial bilinear estimate due to Coifman and Meyer which uses the Carleson measure characterization of BMO functions by Fefferman.

\medskip

Inspired by parametrisation commonly used in meteorology and oceanography, we consider an
alternative approach to the deterministic theory. While it departs from the Rayleigh-B\'enard framework,
it retains the underlying transport mechanism. In the Stokes equations \eqref{eq:UStokes}, we
replace the temperature $T$ by a model $S$ of its fluctuations, given by a noise that is white in
time and coloured in space,
\begin{equation*}
S=\sqrt{\mathcal{Q}_S}\,\dot W,
\end{equation*}
where $\dot W$ is the formal time derivative of a cylindrical Wiener process $W$, i.e.\ a space-time white noise, and $\mathcal{Q}_S$ is a covariance operator,
which encodes the spatial structure of the fluctuations. The velocity field then solves
\[
-\Delta\ub+\nabla p=\Ra\,S\,\bm{e}_{x_3},\qquad\nabla\cdot\ub=0,
\]
in the distributional sense, and is therefore itself a noise of the form\footnote{Formally, applying
the divergence to the Stokes equations and using $\nabla\cdot\ub=0$ gives $\Delta p=\Ra\,\partial_3S$,
and applying the Laplacian gives $-\Delta^2\ub+\nabla\Delta p=\Ra\,\Delta S\,\bm{e}_{x_3}$. Combining the
two,
\[
\ub=-\Ra\,\Delta^{-2}\begin{pmatrix}-\nabla'\partial_3\\ \Delta'\end{pmatrix}S,
\]
where $\nabla'$ and $\Delta'$ denote the horizontal gradient and the horizontal Laplacian, and
$\Delta^{-2}$ is the inverse of the bi-Laplacian with clamped boundary conditions. For the vertical
component, this is the standard reduction $\Delta^2u_3=-\Ra\,\Delta'S$, with $u_3=\partial_3u_3=0$ at the
plates.}
\[
\ub=\sqrt{\mathcal{Q}_\ub}\,\dot W .
\]
On the periodic strip $\D^{{\rm st},L}=\TT^2_L\times(0,1)$ (with $L\geq1$ the period in the horizontal directions), we thus consider the following
(non-dimensionalised) model, in which the buoyancy forcing is parametrised by $S$:
\begin{align}
\partial_tT+(\ub\cdot\nabla)T&=\Delta T,\label{eq:T}\\
\nabla\cdot\ub&=0,\label{eq:div}\\
-\Delta\ub+\nabla p&=\Ra\,S\,\bm{e}_{x_3},\label{eq:u}\\
T(t,\xx',0)&=1,\quad T(t,\xx',1)=0, && t\ge0,\ \xx'\in\TT^2_L,\label{eq:bc T}\\
\ub(t,\xx',0)&=0,\quad \ub(t,\xx',1)=0, && t\ge0,\ \xx'\in\TT^2_L,\,.\label{eq:bc u}
\end{align}
 Our aim is to identify
assumptions on $\mathcal{Q}_S$ which respect the phenomenological structure of convection, and to
understand how they determine the upper bounds on the Nusselt number.

\medskip

The model above is a variant of a classical procedure to study the effect of a turbulent flow on an
advected quantity, namely a decoupling argument consisting of two steps:
\begin{enumerate}
\item[(i)] from the fluid equations, one extracts a synthetic model of the velocity field $\ub$;
\item[(ii)] one inserts this $\ub$ into the equation for the advected quantity.
\end{enumerate}
This strategy has been used in many different contexts, for instance for passive scalars
\cite{kraichnan1994anomalous}, polymers \cite{picardo2023polymers}, magnetic fields
\cite{kazantsev1968enhancement} and temperature in plasmas \cite{vlad1998diffusion}. Most of the
time, the synthetic $\ub$ is chosen on the basis of general properties of the flow; an example in
which it is identified from numerical simulations of the fluid equations can be found in
\cite{cifani2026stochastic}.

For an active advected quantity, such as the temperature in the Rayleigh-B\'enard problem, the feedback of
the temperature on the fluid is an essential part of the mechanism, and a complete decoupling into
the steps (i) and (ii) risks losing it. Our approach keeps this feedback in the simplest possible
form: the synthetic velocity field is governed by the stationary Stokes equations, forced by the
fluctuations of a temperature field. The better the model of these fluctuations, the more
realistic the resulting $\ub$.

This choice has two limitations. The first, and obvious, one is that the loop is open: the model
of the temperature fluctuations is not the one generated by advecting the temperature with $\ub$.
The second is that we only model fluctuations. Since the Stokes equations are applied to the
temperature fluctuations alone, $\ub$ is a model of the turbulent part of the velocity field only,
and the mean temperature, which exists and is central for the Nusselt number, does not enter.
The underlying intuition is that, in the ultimate regime of very large Rayleigh numbers, the flow
is statistically dominated by the fluctuations, and these alone determine the scaling law.

On the other hand, the method has an advantage over choosing the synthetic velocity field a
priori: the boundary layer, that is, the statistical behaviour of $\ub$ near the plates, is not
prescribed, but produced by the Stokes equations from the statistics of the temperature
fluctuations near the plates. In the end, the scaling law depends on the statistical properties of
the turbulent temperature fluctuations near the boundary, and this is will be encoded in our assumptions on
$\mathcal{Q}_S$.

Finally, to avoid a common misunderstanding: the model does not contain any external noise. More precisely, the field $S$
in the Stokes equations is not an external driver or a random perturbation of the system, but a
simplified model of the fluctuations of the internal temperature field, which in the real system
is coupled to the velocity. Instead of computing these fluctuations from the system, we
\emph{choose} a model for them; from it we deduce a model of the velocity, which we then use in the
advection equation to obtain information on the mean temperature.

\medskip

Let $q_S^{{\rm st}, L}$ denote the integral kernel of $\mathcal{Q}_S$, defined by
\begin{equation*}
(\mathcal{Q}_S\phi)(\xx) = \int_{\D^{{\rm st}, L}} q_S^{{\rm st}, L}(\xx,\yy)\phi(\yy){\rm d} \yy.
\end{equation*}
Specifically, our objective is to identify sufficient conditions on $q_S$ to reproduce both the Spiegel–Kraichnan as well as the Malkus scaling. Our main result reads as follows:

\begin{theorem}\label{Main Theorem}
Let $L\geq 1$, $\Ra\gg 1$ and assume that the stochastic temperature covariance $q_S^{{\rm st}, L}$ is symmetric, positive semidefinite, satisfies the boundary conditions
\begin{equation*}
    q_{S}^{{\rm st}, L}(\yy,\zz)=0 \qquad \mbox{ whenever } y_3\in\{0,1\} \mbox{ or } z_3\in\{0,1\}\,,
\end{equation*}
and is horizontally stationary
\begin{equation*}
	q_S^{{\rm st}, L}(\yy+(\vv',0), \zz+(\vv',0)) = q_S^{{\rm st}, L}(\yy,\zz) \qquad \forall \vv'\in \R^2, \yy,\zz\in \DsL\,.
\end{equation*}
\begin{description}
		\item[Case 1] \emph{(Uniform correlation length):} Assume that $q_S^{{\rm st}, L}$ satisfies
			\begin{equation*}
			|q_S^{{\rm st}, L}(\yy,\zz)|\leq m(y_3)m(z_3)\sum_{\kk'\in\ZZ^2}(1+|\yy'-\zz'+L\kk'|)^{-N} \,,\quad N>2,
			\end{equation*}
			with
			\begin{equation*}
			m(r)= r(1-r),
			\end{equation*}
			then
		      \begin{equation*}
				\Nu\leq C\Ra^{\frac 12}\,,
			\end{equation*}
            where the constant is independent of $L$ and $\Ra$.
		\item[Case 2]\emph{(Boundary-layer correlations):} Assume that for some some $0<\eta \leq \frac 12$, $q_S^{{\rm st}, L}$ admits the decomposition
			\begin{equation*}
			     q_S^{{\rm st}, L}=q_S^{\rm{bl},\rm{bottom},\eta}+q_S^{\rm{bulk},\eta}+q_S^{\rm{bl},\rm{top},\eta},
			\end{equation*}
            where $q_S^{\rm{bl},\rm{bottom},\eta}, q_S^{\rm{bulk},\eta}, q_S^{\rm{bl},\rm{top},\eta}$  satisfy the following assumptions:
			\begin{itemize}
            \item $q_S^{\rm{bl},\rm{bottom},\eta}$ is localised in $y_3, z_3 \leq \eta$ and satisfies
            \begin{equation*}
			     \left|q_S^{\rm{bl},\rm{bottom},\eta}(\yy,\zz)\right|\leq y_3z_3\sum_{\kk'\in \ZZ^2}\left(1+\frac{|\yy'-\zz'+L\kk' |}{\sqrt{y_3z_3}}\right)^{-N^{\rm bot}_{\rm{bl}}},\, N^{\rm bot}_{\rm{bl}}>2;
			\end{equation*}
            \item $q_S^{\rm{bl},\rm{top},\eta}$ is localised in $y_3, z_3 \geq 1- \eta$ and satisfies
            \begin{equation*}
			     \left|q_S^{\rm{bl},\rm{top},\eta}(\yy,\zz)\right|\leq (1-y_3)(1-z_3)\sum_{\kk'\in \ZZ^2}\left(1+\frac{|\yy'-\zz'+L\kk'|}{\sqrt{(1-y_3)(1-z_3)}}\right)^{-N^{\rm top}_{\rm{bl}}},\, N^{\rm top}_{\rm{bl}}>2;
			\end{equation*}
            \item $q_S^{\rm{bulk},\eta}$ is localised in the bulk $\eta < y_3,z_3 < 1- \eta$ and satisfies 
			\begin{equation*}
			\left|q_S^{\rm{bulk},\eta}(\yy,\zz)\right|\leq \eta^{2}m_{\rm{bulk}}(y_3)m_{\rm{bulk}}(z_3)\sum_{\kk'\in \ZZ^2}(1+|\yy'-\zz'+ L\kk'|)^{-N_{\rm{bulk}}},\, N_{\rm bulk}>2,
			\end{equation*}
	    	with 
	    	\begin{equation}
	    	m_{\rm{bulk}}(r):= (r(1-r))^{1+\sigma}, \quad \sigma >0.
	    	\end{equation}
        \end{itemize}
	    	Then
	    	\begin{equation*}
				\Nu\leq C \left(\Ra^{\frac 13}+(\Ra \eta)^{\frac 12}\right)\,,
			\end{equation*}
            where the constant is independent of $L, \eta$ and $\Ra$. In particular, 
            \begin{equation*}
            \Nu\leq C\Ra^{\frac 13}\,
            \end{equation*}
            when $\eta\leq \Ra^{-\frac 13}$.
\end{description}
    
\end{theorem}

We call $\eta$ the thickness of the \emph{fluctuation boundary layer}, i.e. the layer in which the
statistics of the temperature fluctuations are shaped by the wall, with an amplitude vanishing
linearly and a correlation length comparable to the distance from the wall. It should not be
confused with the thermal boundary layer, which is not prescribed in the model but is determined
by the heat flux: in our argument, it is the thickness $\delta$ of the averaging layer, with
$\Nu\lesssim 1/\delta$. The Malkus scaling holds precisely when the fluctuation boundary layer is no
thicker than the thermal one, $\eta\lesssim\Ra^{-1/3}$.

In Case 1 of Theorem \ref{Main Theorem}, the bound $\Nu\lesssim\Ra^{1/2}$ arises because the horizontal
correlation structure of the temperature fluctuations does not see the plates. To see this, it is convenient to write the
assumption of Case~1, up to periodisation, as
\[
\left|q_S^{{\rm st}, L}(\yy,\zz)\right|\lesssim m(y_3)\,m(z_3)\left(1+\frac{|\yy'-\zz'|}{\ell_{\rm corr}}\right)^{-N},\quad N>2,
\]
where $\ell_{\rm corr}$ is the horizontal distance over which the fluctuations are correlated. In
Case 1, it holds $\ell_{\rm corr}=1$: the amplitude of the covariance vanishes at the plates, through the
factor $m(y_3)m(z_3)$, but the horizontal correlation length does not see the boundary.
Fluctuations right next to a plate are correlated over the same horizontal distances as those in
the well-mixed bulk, as if the plates do not impose any constraint on the horizontal structure of the
flow, only on its vertical component. This is reminiscent of the fact that, for the classical Rayleigh-B\'enard convection  model, the no-penetration condition alone leads to the bound $\Ra^{1/2}$, and which is at the basis of the so-called direct method \cite{chernyshenko2022relationship}
developed in \cite{otto2011rayleigh,choffrutNobiliOttoUpperBoundsOnNusseltNumberAtFinitePrandtlNumber,seis2015scaling,nobili2023role}. From this point of view, the bound $\Ra^{1/2}$
is obtained ignoring the structure of the boundary layer; to obtain $\Ra^{1/3}$, this structure has to be encoded in the covariance.

\begin{figure}
\centering
\includegraphics[width=\linewidth]{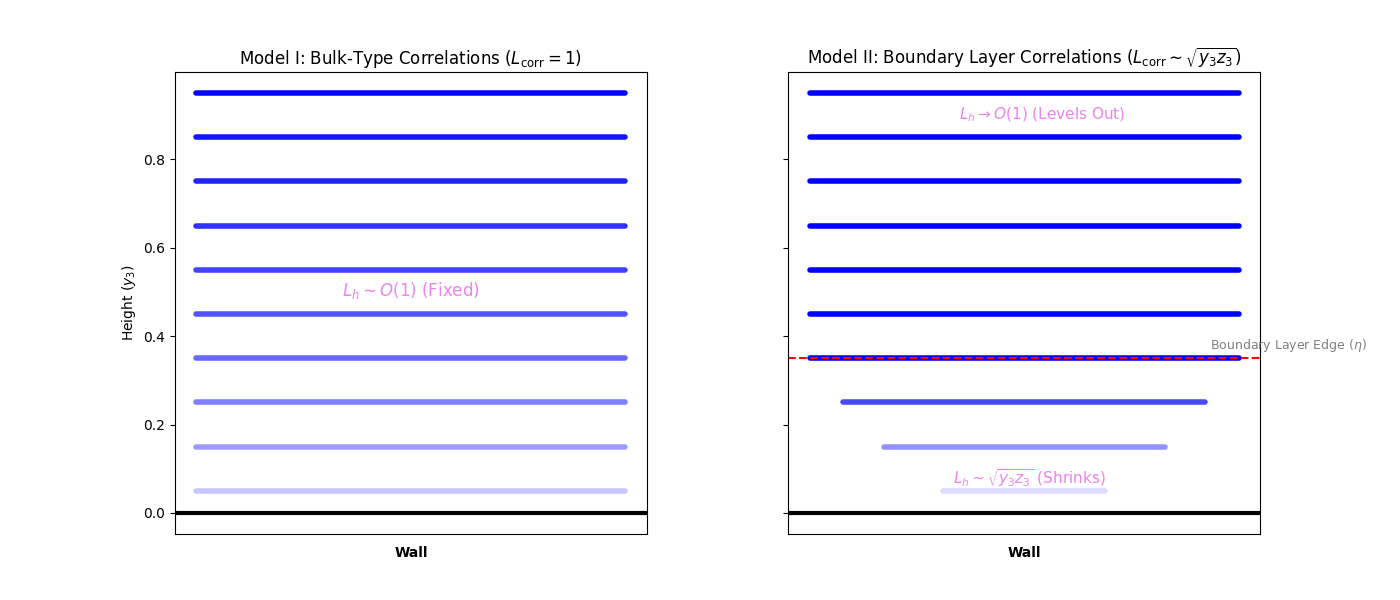}
\caption{Horizontal correlation lengths in the two cases of Theorem~\ref{Main Theorem}. Left (Case~1):
the correlation length is of order one up to the plates. Right (Case~2): in a boundary layer of
thickness $\eta$, the correlation length $\sqrt{y_3z_3}$ shrinks to zero at the plates.}
\label{fig:correlations}
\end{figure}

This is done in Case~2, by treating a boundary layer of thickness $\eta$ and the bulk
separately; see Figure~\ref{fig:correlations}. Within the bottom boundary layer, $y_3,z_3\le\eta$, (and analogously in the top boundary layer), the correlation length is $\ell_{\rm corr}=\sqrt{y_3z_3}$, hence it shrinks to zero at the plates:
fluctuations near a plate are correlated only over horizontal distances comparable to their
distance from the plate, in line with the picture of a no-slip boundary layer. In the bulk, the
correlation length is again of order one, but the covariance carries the prefactor $\eta^2$, in order for
bulk fluctuations to be weak when the boundary layer is thin. The resulting bound,
$\Nu\lesssim\Ra^{1/3}+(\Ra\,\eta)^{1/2}$, shows that the boundary layer determines the heat
transport as soon as $\eta\lesssim\Ra^{-1/3}$, and yields the Malkus scaling $\Ra^{1/3}$; for
$\eta\sim1$, it reduces to the bound $\Ra^{1/2}$ of Case~1.

These conditions may also be compared with the background field method
\cite{DC96}. There, one prescribes the mean temperature profile and controls the
fluctuations through a stability condition, and the bound is of order $1/\delta$ for a profile with
boundary layers of thickness $\delta$, the smallest admissible $\delta$ being determined by
marginal stability. Here, the roles are reversed: we prescribe the statistics of the fluctuations
and leave the mean profile free. Nevertheless, the same structure appears. The conductive
contribution in our proof is of order $1/\delta$, where $\delta$ is the thickness of the layer over
which we average, the critical thickness $\eta\sim\Ra^{-1/3}$ is that of a marginally stable boundary
layer, and the smallness of the bulk fluctuations corresponds to the nearly constant bulk profile of
the optimal backgrounds.

As a technical side note, we remark that the exponent $\sigma>0$ in the bulk profile only serves to avoid logarithmic corrections in the
bound; see Remark \ref{rem:sigma}.

The description obtained here is consistent with the rigorous upper bounds of
\cite{otto2011rayleigh,NO17,ChanMal2025}. We stress, however, that our bounds are not universal: they hold for a
model, and we design the spatial covariance of the stochastic temperature fluctuations so as to
reproduce the two most prominent scaling laws, $\Ra^{1/2}$ and $\Ra^{1/3}$. In this sense, the approach
is close in spirit to the optimal design problems of Doering et al.\ \cite{doering2019optimal} and
Tobasco and Doering \cite{TD2017, doering2019optimal}, which was recently generalized to the stochastic context by the same authors of this paper \cite{FLN2026a}, offering an indirect perspective on the ongoing
debate on the ultimate regime \cite{doering2019,lohse2024ultimate,shishkina2024ultimate}. We finally remark that in \cite{FLN2026c} we used the same novel approach, namely modelling the
temperature fluctuations by a noise that drives the velocity field, for porous medium convection.
There, the stochastic velocity field solves Darcy's law forced by stochastic thermal fluctuations,
and we reproduce the bound $\Nu\lesssim\Ra$ known for the deterministic problem \cite{doering1998bounds}.

\paragraph{Notations}
We will use the following average notation:
\begin{itemize}
	\item long-time average
	\[\overline{(\cdot)} := \limsup_{t_0\rightarrow \infty}\frac{1}{t_0} \int_0^{t_0} (\cdot){\rm d}t;\]
	\item horizontal-space average
	\[\langle (\cdot)\rangle' := \frac{1}{L^2} \int_{\mathbb{T}_L^2} (\cdot) {\rm d} \xx'.\]
\end{itemize}

Our equations are set in the horizontally periodic strip, nevertheless the proof of Proposition \ref{Pr1} is complex and requires us to pass to various extension of this domain. In particular, we will use the following:

\begin{eqnarray*}
    \DsL&=\TT_L^2\times(0,1) \quad &\text{periodic strip},\\
    \DhL&=\TT^2_L\times (0,\infty) \quad &\text{periodic half-space},\\
    \Dis&=\R^2\times(0,1) \quad &\text{infinitely extended horizontal-strip},\\
    \Dh&=\R^2\times(0,\infty) \quad &\text{half-space}.
\end{eqnarray*}

\paragraph{Structure of the paper}
The Nusselt number is controlled by the vertical velocity covariance $q^{{\rm st},L}_{3,3}$ near the plates, therefore the core of the paper consists of a bound on this quantity, uniform in $L$ (Proposition \ref{Pr1}). Since the kernel in the periodic strip is not explicit, we transfer the problem to the half-space: we extend $q^{{\rm st},L}_S$ by zero (Lemma \ref{extension-lemma}) and unfold the periodic covariance into an integral over $\Dh$ with explicit half-space kernel
(Lemma \ref{unfolding-lemma}), which we then estimate (Lemmas \ref{lemma:bulkCov} and \ref{lemma:BL-bulk}). The bounds are transferred back to the strip by a corrector estimate based on energy methods (Section \ref{sec:corrector}), and the bound on the Nusselt number follows by optimising the thickness of the averaging layer (Section \ref{sec:proofMT}).

We think it is helpful to summarize the main steps of the (rather involved) proof of Proposition \ref{Pr1} in a diagram:
\begin{figure}[h]
\centering
\begin{tikzcd}[column sep=4.5em, row sep=3.6em]
\nd{$\DsL$}{periodic strip}
  \arrow[r, "\aname{extension}", "\lref{Lemma \ref{extension-lemma}}"']
&
\nd{$\DhL$}{periodic half-space}
  \arrow[r, "\aname{unfolding}", "\lref{Lemma \ref{unfolding-lemma}}"']
&
\nd{$\Dh$}{half-space}
  \arrow[dl, "\twolab{kernel estimates}{Lemmas \ref{lemma:bulkCov}, \ref{lemma:BL-bulk}}"]
\\
\nd{bound on $q^{{\rm st},L}_{3,3}$}{Proposition \ref{Pr1}}
  \arrow[d, "\twolab{choice of $\delta$}{Section \ref{sec:proofMT}}"']
&
\nd{bound on $q^{{\rm hs},L}_{3,3}$}{}
  \arrow[l, "\aname{corrector}"', "\lref{Lemma \ref{lem:stoch-corr}}"]
&
\\
\nd{$\mathrm{Nu}\lesssim\Ra^{1/2}$, \ $\mathrm{Nu}\lesssim\Ra^{1/3}$}{Theorem \ref{Main Theorem}}
& &
\end{tikzcd}
\label{fig:roadmap}
\end{figure}

The upper row reduces the problem, posed in the periodic strip,
to the half-space, where the kernel is explicit; the bounds on the vertical velocity covariance
obtained there are then carried back to the strip.

\section{Stochastic framework}\label{sec:stoch-frame}
In this section we will motivate our notion of temperature fluctuations
$$S=\sqrt{\QQ_S}\dot{W}$$
and rigorously implement the stochastic framework we will be working in throughout this manuscript. Furthermore, we would like to stress that the computations in the following sections are rigorously implemented only when considered on the time-integrated level of Brownian motion, instead of its distributional time-derivative. Yet we note that the evolution equations for $\ub$ are elliptic (and $S$ is the only source of time dependence), in which the formal setting of $\dot{W}$ can be easily justified and we remain in that set-up.

\subsection{A stochastic model for temperature fluctuations}\label{sec:TempFluc}

We model our intuition on the temperature fluctuations by imposing that $S$ shall be a centred Gaussian field, statistically stationary in time, statistically invariant by horizontal shifts, with a vertical damping at $x_{3}=0$ and $x_3=1$. More precisely, $S$ shall be described by the form
\[ S(t,\xx) = f(x_3)S_0(t,\xx)\]
where $S_0(t,\xx)$ is a centred, completely homogeneous Gaussian field, namely defined for \textit{all} $\xx \in\mathbb{R}^{3}$ and invariant by \textit{all} shifts, in time and space, and $f(x_3)$ is a damping factor, described below by condition (\ref{condition on f}).\\
\noindent
Let us discuss the derivation of our modelling assumptions in detail. We proceed in two steps:
\begin{enumerate}
	\item First assume a short time correlation $\mathfrak{t}$ for $S_0$ with exponential decay (in order to work with real-valued, not distributional, processes of finite variance): denote $C$ the covariance kernel of $S_0$, then a model for
$S_0$ is a stationary solution of the Ornstein-Uhlenbeck equation
\begin{equation}\label{eq:OU}
	{\rm d}S_0(t,\xx)  = -\frac{1}{\mathfrak{t}}S_0(t,\xx) {\rm d}t + \frac{\sqrt{2}}{\sqrt{\mathfrak{t}}}{\rm d}B(t,\xx)
\end{equation}
where $B$ is a Brownian motion in time with space correlation $C(\xx,\yy)$; the explicit form is
\begin{equation}\label{eq:explicitOUformula}
	S_0(t,\xx) = \int_{-\infty}^{t} e^{-\frac{t-s}{\mathfrak{t}}}\frac{\sqrt{2}}{\sqrt{\mathfrak{t}}}{\rm d}B(s,\xx)
\end{equation}
where (by abuse of notation) $B$ is a two-sided extension in time of the Brownian motion. Note that as stated before, it is rather natural to assume homogeneity in the sense of $C(\xx,\yy) = C(\xx - \yy)$.\\
\noindent
The function $f(x_3)$ plays a crucial role: the full temperature $T$, given by the sum of an average term $\overline{T}$ and the fluctuations $S$, must satisfy the deterministic boundary conditions
\[T(t,\xx',0)  =1 \quad \text{and} \quad T(t,\xx',1)=0,\]
and the average term $\overline{T}(t,\xx)$ therefore has to satisfy the same condition, $\overline{T}(t,\xx',0) = 1, \overline{T}(t,\xx',1) = 0$. Consequently, for the fluctuations, we have the deterministic condition
\[S(t,\xx',0) = S(t,\xx',1) =0.\]
The decay to zero of the process $x_{3}\mapsto S(t,\xx',x_3)$ as $x_{3}\rightarrow0$ or $x_{3}\rightarrow1$ is a complicated issue with fundamental consequences on the full model. The average profile $\overline{T}$, as we see a posteriori from the Nusselt number computation, approaches linearly the limit value $\overline{T}(t,\xx',0) =1$, namely
\[ \lim_{x_3\rightarrow0}\frac{\overline{T}(t,\xx',x_3)-1}{x_3} \in (-\infty,0).\]
One can argue similarly in the case $x_3 \rightarrow 1$. Therefore the most natural assumption is that also the fluctuations decay linearly, based on the general vague practical rule that ``the fluctuations are roughly proportional to the trend". Based on these heuristic arguments, we assume that $f$ is a continuous real function such that
\begin{equation}\label{condition on f}
	0<f(x_3)\leq1, x_3>0\text{, non decreasing, and}\lim_{x_3\rightarrow0}\frac{f(x_3)}{x_3} \geq 0, \lim_{x_3\rightarrow1}\frac{f(x_3)}{x_3} \geq 0.
\end{equation}
	\item Idealisation to white-in-time noise ($\mathfrak{t} \rightarrow 0$): dealing with positive time correlation renders the mathematical analysis rather difficult. The quantitative step from short to delta correlation is based on the following approximation
\[S_0(t,\xx) =\sqrt{2\mathfrak{t}}\frac{{\rm d}B(t,\xx)}{{\rm d}t}.\]
It is heuristically immediately justified by setting the right-hand-side of equation \eqref{eq:OU} equal to zero; otherwise, one may invoke the following simple inequality
\[\mathbb{E}\left[  \left\vert \frac{1}{\sqrt{2\mathfrak{t}}}\int_{0}^{t}S_0(s,\xx) {\rm d}s - B(t,\xx)\right\vert ^{2}\right]\leq \mathfrak{t},\]
easy to prove by elementary manipulations of the formula \eqref{eq:explicitOUformula}. Therefore we arrive at the model
\[S(t,\xx) = \sqrt{2\mathfrak{t}}f(x_3) \frac{{\rm d}B(t,\xx)}{{\rm d}t}.\]
The stochastic process $\sqrt{2\mathfrak{t}}f(x_3)B(t,\xx)$ has time-one covariance equal to
\begin{equation}\label{eq:covMotivation}
    \mathbb{E}\left[  \sqrt{2\mathfrak{t}}f(x_3)B(1,\xx)\sqrt{2\mathfrak{t}}f(x_3)B(1,\yy)\right]  =2\mathfrak{t} f(x_3)f(y_3) C(\xx,\yy).
\end{equation}
\end{enumerate}
Therefore the final model is
\begin{equation}\label{eq:Sform}
    S = \sqrt{\mathcal{Q}_S}\dot{W}
\end{equation}
where $\dot{W}$ is a space-time white noise (with covariance given by the distribution $\delta(\xx-\yy)$), and $\mathcal{Q}_S$ a spatial covariance operator in particular encoding the behaviour towards the boundary as motivated in \eqref{eq:covMotivation}. We will assume that $\mathcal{Q}_S$ is a Hilbert-Schmidt operator, thus denoting $q_S^{{\rm st}, L}$ the integral kernel of $\mathcal{Q}_S$, i.e.
\[(\mathcal{Q}_S \phi)(\xx) = \int_{\D^{\rm st, L}} q_S^{{\rm st}, L}(\xx,\yy)\phi(\yy) {\rm d} \yy,\]
then $q_S^{{\rm st}, L}$ satisfies the following immediate properties:
\begin{itemize}
	\item symmetric, and positive semi-definite, i.e.
		\[\int_{\D^{{\rm st},L}} \int_{\D^{{\rm st},L}} q_S^{{\rm st}, L}(\xx,\yy)\phi(\xx)\phi(\yy) {\rm d} \xx {\rm d}\yy \geq \|\phi\|^2_{L^2(\D^{{\rm st},L})} \quad \forall \phi \in L^2\left(\D^{{\rm st},L}\right);\]
	\item $q_S^{{\rm st}, L} \in L^2(\D^{{\rm st}, L} \times \D^{{\rm st},L})$;
	\item $\widetilde{q}^{{\rm st}, L}_S(\cdot) := q_S^{{\rm st}, L}(\cdot,\cdot) \in L^1(\D^{{\rm st}, L})$.
\end{itemize}
In view of the above discussion, it holds
\[q_S^{{\rm st}, L}(\xx,\yy) = 2\tau f(x_3) f(y_3) C(\xx,\yy)\]
under condition (\ref{condition on f}) on $f$, and in particular
\begin{equation*}
    q_S^{{\rm st}, L}(\xx,\yy) = 0 \quad \text{ whenever } x_3 \in \{0,1\} \text{ or } y_3 \in \{0,1\}.
\end{equation*}

\subsection{Notion of solutions}\label{sec:solutions}
As discussed in the introduction, a noise \eqref{eq:Sform} together with the Stokes equations gives rise to a velocity field of the form
\[\ub = \sqrt{\mathcal{Q}_\ub}\dot{W}.\]
Analogously to \cite[Section 2]{FLN2026a}, this shall denote a cylindrical Brownian motion on a probability space $(\Omega, \mathcal{F}, \mathbb{P})$ with respect to a divergence-free complete orthonormal system, with spatial covariance operator $\mathcal{Q}_\ub$, i.e.
\begin{equation*}
    {\rm d} \ub = \sqrt{\mathcal{Q}_\ub}(\xx){\rm d}W_t = \sum_{k \in \mathbb{N}} \lambda_k e_k(\xx) {\rm d}\beta_t^k
\end{equation*}
where $(\lambda_k) \in \ell^2$, $(e_k)_k$ is a complete orthonormal system of $L_\sigma^2\left(\D^{{\rm st}, L}\right) =$\\ $ \overline{\{f \in C_c^\infty: \nabla \cdot f =0\}}^{L^2(\D^{{\rm st}, L}; \mathbb{R}^3)}$ with $e_k \in C\left(\overline{\D^{{\rm st}, L}}; \mathbb{R}^3\right)$ and $ e_k(x_3=0,1) = 0$, while $\beta^k$ are independent one-dimensional standard Brownian motions. The kernel $Q_\ub$ of $\mathcal{Q}_\ub$, defined by
\[Q_\ub(\xx,\yy) = \sum_k \lambda_k^2 e_k(\xx) \otimes e_k(\yy),\]
therefore is symmetric, positive semi-definite, $Q_\ub \in L^2\left(\D^{{\rm st}, L} \times \D^{{\rm st}, L}\right)$ and satisfies the boundary conditions
\[Q_\ub(\xx,\yy) = 0 \text{ whenever } x_3 \in \{0,1\} \text{ or }y_3 \in \{0,1\}.\]
Such a velocity field gives rise to the following stochastic advection-diffusion equation for the temperature field
\begin{equation*}
	{\rm d} T + \left(\left(\sqrt{\mathcal{Q}_\ub}\circ {\rm d}W_t\right)\cdot \nabla \right) T = \Delta T {\rm d}t
\end{equation*}
in Stratonovich form, or equivalently for $\tilde{Q}_\ub(\cdot) := Q_\ub(\cdot,\cdot)$
\begin{equation}\label{eq:Ito}
	{\rm d} T + \left(\left(\sqrt{\mathcal{Q}_\ub} {\rm d}W_t\right)\cdot \nabla \right) T = {\rm div}\left(\left({\rm Id} +\frac{1}{2} \tilde{Q}_\ub \right) \nabla T\right){\rm d}t
\end{equation}
in It\^{o} form (e.g. c.f. Section 3.4 in \cite{CG24}). Let $\tilde{Q}_\ub(\xx) = [\tilde{q}_{i,j}(\xx)]_{i,j=1,\dots,3}$ and assume that
\begin{equation}\label{eq:regTildeQ}
    \tilde{Q}_\ub \in W^{1,\infty}\left(\D^{{\rm st}, L}\right).
\end{equation}
For the following spaces
\begin{align*}
        H &:= H_0^1\left(\D^{{\rm st},L}\right), \text{ with inner product } \langle \cdot, \cdot \rangle_H\\
        V &:= H^2\left(\D^{{\rm st}, L}\right) \cap H_0^1\left(\D^{{\rm st}, L}\right),
    \end{align*}
one can deduce the existence of a self-adjoint operator $\mathcal{A}$, such that $- \mathcal{A}$ generates a strongly continuous semi-group $e^{t\mathcal{A}}, t \geq 0$, on $H$ (for details c.f. \cite[Section 2]{FLN2026a}).
\begin{definition}[Definition of probabilistic solution]\label{def:solT}
	Let $\tau \geq 0$ and let $\mathfrak{T}= \mathfrak{T}(x_3)$ be a smooth deterministic function satisfying $\mathfrak{T}(0)=1, \mathfrak{T}(1)=0$. Given $T_0$ such that $\mathcal{T}_0 := T_0 - \mathfrak{T} \in H$, we say that $T$ is a \emph{probabilistically strong, mild solution} to \eqref{eq:Ito} with initial condition $T_0$, if $\mathcal{T}:= T - \mathfrak{T}$ 
    is progressively measurable with paths of class $L^\infty(0,\tau;H) \cap L^2(0,\tau;V)$, and satisfies
	\begin{equation}\label{eq:mildSol}
		\mathcal{T}(t) = e^{t\mathcal{A}}\mathcal{T}_0 - \int_0^t e^{(t-s)\mathcal{A}}\left(\sqrt{\mathcal{Q}_\ub} {\rm d} W_s\right)\cdot \nabla \mathcal{T}(s) - \int_0^t e^{(t-s)\mathcal{A}}\left(\sqrt{\mathcal{Q}_\ub} {\rm d} W_s\right)\cdot \mathfrak{T}'\bm{e}_{x_3}
	\end{equation}
	with homogeneous boundary conditions, as well as the energy bound
	\begin{equation}\label{eq:energyInequ}
		\mathbb{E}\left[\sup_{s \leq t}\|\mathcal{T}(s)\|^2_{L^2(\D^{{\rm st},L})}\right] + \mathbb{E}\left[\int_0^t \left\|\nabla \mathcal{T}(s)\right\|_{L^2(\D^{{\rm st},L})}^2{\rm d}s\right] < \infty.
	\end{equation}
\end{definition}
\noindent
The existence of such solutions is rather classical and we refer e.g. to Chapter 3 in \cite{FL2023}, and references therein. Furthermore we observe that the expected value of the It\^{o}-integrals in \eqref{eq:mildSol} vanishes due to square integrability coming from \eqref{eq:energyInequ} thus following from our assumptions on the vector fields $(e_k)_k$ and using the regularity assumption \eqref{eq:regTildeQ} we deduce that $E:= \mathbb{E}[T]$ is a strong solution to
	\begin{align}
		&\partial_t E =\nabla \cdot\left(\left({\rm Id} + \frac{1}{2}\tilde{Q}_\ub\right)\nabla E\right),\label{eq:EvolExp}\\
        &E(t,\xx',0) =1, \quad E(t,\xx',1) = 0,\label{eq:EvolExpBC}
	\end{align}
satisfying $E \in L^\infty(0,\tau;H^1)\cap L^2(0,\tau;H^2), \partial_t E \in L^2(0,\tau;L^2)$.\\
Similar to \cite[Lemma 3.1]{FLN2026a}, \eqref{eq:EvolExp}-\eqref{eq:EvolExpBC} finally allow us to define the Nusselt number
\begin{equation}\label{eq:stochasticNu}
    \Nu := -\left.\overline{\left\langle \partial_3 \mathbb{E}[T]\right\rangle' }\right|_{x_3 =0}.
\end{equation}

\section{Preliminaries}

\subsection{A general upper bound for the Nusselt number}
\begin{lemma}\label{lem:Nusselt}
	Let $\delta \in (0,1]$ be arbitrary, then the Nusselt number \eqref{eq:stochasticNu} satisfies
    \begin{equation}\label{eq:Nugeneral}
        \Nu\lesssim \frac1\delta+\frac1{\delta^2}\sum_{j=1}^3\int_0^\delta
\left\langle|\widetilde q_{3,j}|^2\right\rangle'(x_3)\,{\rm d}x_3.
    \end{equation}
    If $\widetilde{q}_{3,j} = \widetilde{q}_{3,j}(x_3), j=1,\dots 3$, then this can be improved to
	\begin{equation}\label{Nu-apriori}
		\Nu \lesssim \frac{1}{\delta} + \frac{1}{\delta^2}\int_0^{\delta} |\widetilde{q}_{3,3}(x_3)|^2{\rm d}x_3.
	\end{equation}
\end{lemma}
The bound \eqref{eq:Nugeneral} is rather general and does not require more assumptions than the ones in Section \ref{sec:solutions}. Its proof can be found in \cite[Theorem 4.1]{FLN2026a}. By Lemma \ref{ap:lemma1/2} below, our framework here satisfies the more specific assumption of \eqref{Nu-apriori}, hence we include its proof in Appendix \ref{app:nusselt}.

\subsection{Kernels}

Let $K(\xx,\yy)=\{K_{i,k}(\xx,\yy)\}_{i,k=1,\cdots,3}$ denote the $3\times 3$ matrix of kernels associated to \eqref{eq:u} and \eqref{eq:bc u}. In the deterministic case, for any vector-valued forcing $\phi=(\phi_1,\phi_2,\phi_3)$ we have 
\[u_{i}(t,\xx)=\sum_{k=1}^{3}\int_{\D^{{\rm st},L}} K_{i,k}(\xx,\yy)\phi_k(t,\yy)\,{\rm d}\yy\,.\]
In our case $\phi(t,\yy)=(0,0,S(t,\yy))$, hence formally it holds,
$$u_i(t,\xx)=\Ra\int_{\D^{{\rm st},L}} K_{i,3}(\xx,\yy)S(t,\yy)\, {\rm d}\yy\,.$$
Rigorously we will consider the entries of the covariance kernel $Q_{\ub}$, given by
\begin{equation*}
	q_{i,j}(\xx,\yy)=\Ra^2 \int_{\D^{{\rm st},L}}\int_{\D^{{\rm st},L}} K_{i,3}(\xx,\xib)q_S^{{\rm st},L}(\xib,\etab)K_{j,3}(\etab,\yy)\, {\rm d}\xib\, {\rm d}\etab.
\end{equation*} 

\subsection{Horizontal stationarity}
Horizontal shifts $\xib\mapsto\xib+(\vv',0)$, $\vv'\in\mathbb{R}^2$, are understood modulo $L$ in the horizontal variables.
Observe that the kernel $K$ is translation invariant in the horizontal directions, i.e.
\begin{equation}\label{folk}
    K(\xib+(\vv',0), \etab+(\vv',0)) = K(\xib,\etab)
	\qquad\quad \forall \vv'\in \R^2, \,\forall \xib, \etab \in \D^{{\rm st},L}.
\end{equation}
We consider now the translated coefficients
\begin{align*}
	&q_{i,j}(\xx+(\vv',0),\yy+(\vv',0))\\
	&=\Ra^2 \int_{\D^{{\rm st},L}}\int_{\D^{{\rm st},L}} K_{i,3}(\xx+(\vv',0),\xib)q_S^{{\rm st},L}(\xib,\etab)K_{j,3}(\yy+(\vv',0), \etab)\, {\rm d}\xib\, {\rm d}\etab\,.
\end{align*}
Assuming that also $q_S^{{\rm st},L}$ is horizontally stationary in the sense that
\begin{equation}\label{qS:HorStat}
    q_S^{{\rm st},L}(\xib+(\vv',0), \etab+(\vv',0)) = q_S^{{\rm st},L}(\xib,\etab)
	\qquad\quad \forall \vv'\in \R^2, \,\forall \xib, \etab \in \D^{{\rm st},L},
\end{equation}
we find that by a change of variables $\xib=\bar{\xib}+(\vv',0)$, $\etab=\bar{\etab}+(\vv',0)$, and the property \eqref{folk} we deduce that $q_{i,j}$ are horizontally stationary. In particular, setting $\xx=\yy$ it holds
 \begin{equation*}
 q_{i,j}(\xx+(\vv',0),\xx+(\vv',0))=q_{i,j}(\xx,\xx).
 \end{equation*}
Hence, setting $\vv'=-\xx'$ we have 
$$q_{i,j}((\xx',x_3),(\xx',x_3) )=q_{i,j}((\bm{0}',x_3), (\bm{0}',x_3)).$$
This way we deduce the following result:
\begin{lemma}[Horizontal stationarity]\label{ap:lemma1/2}
	Assume that $q_S^{{\rm st},L}$ satisfies \eqref{qS:HorStat}, then the coefficients $q_{i,j}, i,j=1,\dots,3,$ are invariant under horizontal translations 
	\begin{equation*}
		q_{i,j}(\xx+(\vv',0),\yy+(\vv',0))=q_{i,j}(\xx,\yy) \qquad \forall \vv'\in \R^2, \,\forall \xx,\yy\in \D^{{\rm st},L}\,.
	\end{equation*}
	Consequently, the diagonal coefficients 
	\begin{equation*}
		\widetilde{q}_{i,j}(\xx):=q_{i,j}(\xx,\xx)
	\end{equation*}
	do not depend on the horizontal variable $\xx'$.
\end{lemma}
In particular, we will denote
\[\widetilde{q}_{i,j}(x_3) : = \widetilde{q}_{i,j}(\bm{0}',x_3) = q_{i,j}((\bm{0}',x_3),(\bm{0}',x_3)).\]

Note that the same argument applies verbatim on the periodic half-space $\D^{{\rm hs},L}$, with $K$ and
$q_S^{{\rm st},L}$ replaced by $K^{{\rm hs},L}$ and $q^{{\rm hs},L}_S$; in particular, $q^{{\rm hs},L}_{3,3}(\xx,\xx)$ depends only on $x_3$. Therefore from Lemma \ref{lem:Nusselt}, in particular \eqref{Nu-apriori}, the upper bound of the Nusselt number will be determined by the behaviour of $\widetilde{q}_{3,3}$.
\subsection{Green's functions on the positive half-space}\label{preliminaries}
In this section, we consider the setting defined on the positive half-space $\HS = \mathbb{R}^2\times(0,\infty)$, letting $S$ be a random forcing on $\HS$ defined via space-time white noise similar to before. By abuse of previous notation, let $K(\xx,\yy)=\{K_{i,k}(\xx,\yy)\}_{i,k=1,\cdots,3}$ denote the $3\times 3$ matrix of kernels associated to \eqref{eq:u} and \eqref{eq:bc u} on $\HS$. In particular, we need to analyse the kernels $K_{i,3}$ with $i=1,2,3$ in detail.
Let us start with $i=3$: it is easy to show that the vertical velocity $u_3 = \ub \cdot \boldsymbol{e}_{x_3}$ satisfies the 4th–order elliptic boundary value problem
\begin{equation}\label{fourth-order-elliptic}
	\begin{cases}
		-\Delta^2 u_3 = \Ra \,\Delta' S & \text{in } \HS \\
		u_3 = \partial_3 u_3 = 0 & \text{at } x_3=0
	\end{cases}
\end{equation}
where $\Delta'$ is the horizontal Laplacian $\Delta'=\partial_1^2+\partial_2^2$.
The knowledge of the kernel $K_{3,3}$ associated to the operator in \eqref{fourth-order-elliptic} is now crucial in order to study the properties of the covariance function of the vertical velocity component. Let  $G(\xib,\etab)$ solve
\begin{equation*}
    \begin{cases}
	-\Delta_{\etab}^2 G(\xib,\cdot) = \delta_\xib & \text{in } \HS \\
	G(\xib,\cdot)=0 & \text{on } \partial\HS \\
	\partial_{\eta_3} G(\xib,\cdot) = 0 & \text{on } \partial\HS,
\end{cases}
\end{equation*}
with $G(\xib,\etab)$ decaying for $\etab_3\rightarrow \infty$. This problem has been extensively studied and the exact expression for $G$ is known (e.g. refer to \cite{ChanMal2025}):
\begin{equation}\label{green-kernel}
	G(\xib,\etab) = -\frac{1}{16\pi}\frac{(|\xib^*-\etab|-|\xib-\etab|)^2}{|\xib^*-\etab|}\,
\end{equation}
where $\xib^{\star} = (\xi_1,\xi_2,-\xi_3).$ Then using that
\begin{align*}
    G(\xib,\etab) = 0 \quad &\text{ at }\eta_3 = 0,\\
    G(\xib,\etab) \rightarrow \infty \quad &\text{ as }\eta_3 \rightarrow 0,\\
    G(\xib,\eta_i), \partial_{\eta_i}G(\xib,\eta_i) \rightarrow 0 \quad &\text{ as }\eta_i \rightarrow \pm \infty,\, i=1,2,
\end{align*}
and that $S=0$ at $\eta_3 =0$, we obtain that the solution of \eqref{fourth-order-elliptic} is given by
\[\int_\HS G(\xib,\etab)\Delta_{\etab}'S(t,\etab)\,{\rm d}\etab=\int_\HS \Delta_{\etab}'G(\xib,\etab)S(t,\etab)\,{\rm d}\etab\,.\]
This representation formula is justified in the weak sense: indeed, we localise in the horizontal variables and define
	$$G_R(\xib,\etab):=\chi_R(\etab')\,G(\xib,\etab),$$
	with $\chi_R\in C_c^\infty(\R^2)$, $\chi_R\equiv 1$ on $\{|\etab'|\le R\}$. For each $R$, the function $G_R(\xib,\cdot)$ is smooth and compactly supported in the horizontal variables, such that
	$$\langle \Delta_{\etab}' S, G_R(\xib,\cdot)\rangle=
	\langle S, \Delta_{\etab}' G_R(\xib,\cdot)\rangle$$
	holds by definition of distributional derivatives. Expanding the above yields
	$$\Delta_{\etab}' G_R=\chi_R \Delta_{\etab}' G
	+2\nabla\chi_R\cdot \nabla_{\etab}'G+(\Delta_{\etab}'\chi_R)G.$$
	Since the commutator terms are supported in the annulus $\{R\lesssim |\etab'|\lesssim 2R\}$, which for fixed $\xib$ lies in the far-field regime as $R\to\infty$, passing to the limit is justified.

One can compute that
\begin{equation*}
	\Delta_{\etab}' G(\xib,\etab) =-
	 \frac{1}{8\pi}\left(\frac{1}{|\xib^*-\etab|} - \frac{1}{|\xib-\etab|^3} - 
	 \frac{(\xi_3-\eta_3)^2}{|\xib-\etab|^3} + \frac{\xi_3^2+\eta_3^2}{|\xib^*-\etab|^3} + 
	 6\frac{\xi_3\eta_3(\xi_3+\eta_3)^2}{|\xib^*-\etab|^5}\,\right)
\end{equation*}
and we define 
\begin{equation}\label{eq:kernelK33}
	K_{3,3}(\xib,\etab) := \Delta_{\etab}' G(\xib,\etab)\,.
\end{equation}
Observe that $K_{3,3}$ enjoys the following three important properties:
\begin{align*}
	K_{3,3}(\xib,\etab) &= K_{3,3}(\etab,\xib)\\
	K_{3,3}(r\xib,r\etab) &= r^{-1}K(\xib,\etab), \qquad \forall r>0,\\
	K_{3,3}(\xib+(\vv',0), \etab+(\vv',0)) &= K_{3,3}(\xib,\etab)
	\qquad\quad \forall \vv'\in \R^2.
\end{align*}

The next lemma contains fundamental near-field/far-field bounds for $K_{3,3}$:
\begin{lemma}[Near-/Far-field kernel bounds]\label{kernel_decay}
	For $\xib\neq \etab$, $\xi_3, \eta_3>0$, the kernel $K_{3,3}$ satisfies:
	\begin{equation}\label{eq:near-field}
		|K_{3,3}(\xib,\etab)|\lesssim \frac{C}{|\xib-\etab|} \qquad \mbox{ for  } |\xib-\etab|\leq 8\min\{\xi_3,\eta_3\}
	\end{equation}
	and 
	\begin{equation}\label{eq:far field-new}
		|K_{3,3}(\xib,\etab)|\lesssim \frac{\min\{\xi_3,\eta_3\}^2}{|\xib-\etab|^3}\qquad \mbox{ for  }  |\xib-\etab|>8\min\{\xi_3,\eta_3\}.
	\end{equation}
\end{lemma}
The proof of this lemma is standard and contained in the appendix (c.f. Lemma \ref{appendix:KEstimates}).\\

For completeness, we determine the structure of the kernel $K_{j,3},\, j=1,2$: the incompressibility condition provides us with
\[\nabla' \cdot \ub' = - \partial_3 u_3\]
where $\ub' = (u_1,u_2)$ and $\nabla' = (\partial_1,\partial_2)^T$, thus in horizontal Fourier variables we obtain 
\begin{equation}\label{id}
	{\rm i}\kk'\cdot \widehat{\ub'}=-\partial_3 \widehat{u_3}
\end{equation}
where ${\rm i}$ denotes the imaginary unit and $\kk' \in \ZZ^2$. Using the decomposition
\begin{equation}\label{eq:FourierHor}
    \widehat{\ub'} = \alpha(\kk',x_3)\kk' + \vv^\perp, \quad \vv^\perp \cdot \kk' =0,
\end{equation}
then inserting into \eqref{id} we obtain
\[\alpha(\kk',x_3) = {\rm i} \frac{\partial_3 \widehat{u_3}}{|\kk'|^2}.\]
Using the Stokes equation for the horizontal part, we obtain in horizontal Fourier variables
\[(|\kk'|^2 - \partial_3^2)\widehat{\ub'} + {\rm i}\kk'\widehat{p} =0\]
and using \eqref{eq:FourierHor}, we obtain
\[(|\kk'|^2 - \partial_3^2)\alpha \kk' + (|\kk'|^2 - \partial_3^2) \vv^\perp + {\rm i}\kk'\widehat{p} =0.\]
In particular, this implies
\[(|\kk'|^2 - \partial_3^2) \vv^\perp \cdot \kk'^\perp = 0 \quad \forall \kk'^\perp \neq 0.\]
Observe that the solution to
\begin{align*}
    (|\kk'|^2 - \partial_3^2) \vv^\perp &= 0,\\
    \vv^\perp(\kk',0) &= 0
\end{align*}
(together with decay assumptions on $\D$) is in fact $\vv^\perp = 0$. This gives in total
\[\widehat{u_j}(\kk',x_3) = \frac{{\rm i} k_j}{|\kk'|^2}\partial_3 \widehat{u_3}(\kk',x_3),\]
which in physical variables gives
\[u_j = \partial_j(\Delta_h)^{-1}\partial_3 u_3.\]
In particular this implies
\[K_{j,3}(\xib,\etab) = \partial_{\eta_j} \partial_{\eta_3} G(\xib,\etab)= \partial_{\eta_j}\partial_{\eta_3} (\Delta_{\etab'})^{-1}K_{3,3}(\xib,\etab).\]

\section{Proof of Theorem \ref{Main Theorem}}
The following proposition provides the core estimates used to prove Theorem \ref{Main Theorem}. Its proof will be tackled in several subsections in which we show that the relevant estimates can be obtained when transferring to a suitable setting on the positive half-space.
\begin{proposition}[Covariance bounds]\label{Pr1}
Recall that the wall-normal velocity covariance be defined as
\[q_{3,3}^{{\rm st}, L}(\xx,\xx)=\Ra^2\int_{\DsL} \int_{\DsL} K^{{\rm st}, L}_{3,3}(\xx,\yy)\,q_S^{{\rm st}, L}(\yy,\zz)\,K^{{\rm st}, L}_{3,3}(\xx,\zz)\,{\rm d}\yy\,{\rm d}\zz,\]
where $K^{{\rm st}, L}_{3,3}$ is the wall-normal component of the Stokes kernel in the periodic strip \eqref{eq:kernelK33}. 
The stochastic temperature covariance $q_S^{{\rm st}, L}$ is assumed to be positive semidefinite, satisfy the boundary conditions
\begin{equation*}
    q_{S}^{{\rm st}, L}(\yy,\zz)=0 \qquad \mbox{ whenever } y_3\in\{0,1\} \mbox{ or } z_3\in\{0,1\}\,,
\end{equation*}
and to be horizontally stationary
\begin{equation*}
	q_S^{{\rm st}, L}(\yy+(\vv',0), \zz+(\vv',0)) = q_S^{{\rm st}, L}(\yy,\zz) \qquad \forall \vv'\in \R^2, \yy,\zz\in \DsL\,.
\end{equation*}
Let $\beta > 0$.
\begin{description}
		\item[Case 1]\emph{(Bulk-dominated):} If the temperature covariance kernel satisfies
			\begin{equation*}
			\left|q_S^{{\rm st}, L}(\yy,\zz)\right|\leq m(y_3)m(z_3)\sum_{\kk'\in\ZZ^2}(1+|\yy'-\zz'+L\kk'|)^{-N} \,,\qquad N>2,
			\end{equation*}
			with
			\begin{equation*}
			m(r)= (r(1-r))^{\beta},
			\end{equation*}
			then, for any fixed $0<x_3\leq 1$, we have
			\begin{equation}\label{q33-bulk}
				q_{3,3}^{{\rm st}, L}(\xx,\xx)\lesssim\Ra^2 x_3^{4}.
			\end{equation}
		\item[Case 2]\emph{(Boundary-layer-dominated):} If instead for $\eta \leq \frac{1}{2}$
			\begin{equation*}
			     q_S^{{\rm st}, L}=q_S^{\rm{bl},\rm{bottom},\eta}+q_S^{\rm{bulk},\eta}+q_S^{\rm{bl},\rm{top},\eta},
			\end{equation*}
            where $q_S^{\rm{bl},\rm{bottom},\eta}, q_S^{\rm{bulk},\eta}, q_S^{\rm{bl},\rm{top},\eta}$ are symmetric and horizontally stationary, the following assumptions hold
			\begin{itemize}
            \item $q_S^{\rm{bl},\rm{bottom},\eta}$ is localised in $y_3, z_3 \leq \eta$ and satisfies
            \begin{equation}
			     \left|q_S^{\rm{bl},\rm{bottom},\eta}(\yy,\zz)\right|\leq (y_3z_3)^\beta\sum_{\kk'\in \ZZ^2}\left(1+\frac{|\yy'-\zz'+L\kk'|}{\sqrt{y_3z_3}}\right)^{-N^{\rm bot}_{\rm{bl}}},\quad N^{\rm bot}_{\rm{bl}}>2;
			\end{equation}
            \item $q_S^{\rm{bl},\rm{top},\eta}$ is localised in $y_3, z_3 \geq 1- \eta$ and satisfies
            \begin{equation}
			     \left|q_S^{\rm{bl},\rm{top},\eta}(\yy,\zz)\right|\leq ((1-y_3)(1-z_3))^\beta\sum_{\kk'\in \ZZ^2}\left(1+\frac{|\yy'-\zz'+L\kk'|}{\sqrt{(1-y_3)(1-z_3)}}\right)^{-N^{\rm top}_{\rm{bl}}},\, N^{\rm top}_{\rm{bl}}>2;
			\end{equation}
            \item $q_S^{\rm{bulk},\eta}$ is localised in the bulk $\eta < y_3,z_3 < 1- \eta$ and satisfies 
			\begin{equation}
			\left|q_S^{\rm{bulk},\eta}(\yy,\zz)\right|\leq \eta^{2\beta}m_{\rm{bulk}}(y_3)m_{\rm{bulk}}(z_3)\sum_{\kk'\in \ZZ^2}(1+|\yy'-\zz'+ L\kk'|)^{-N_{\rm{bulk}}},\, N_{\rm bulk}>2\,,
			\end{equation}
	    	with 
	    	\begin{equation}
	    	m_{\rm{bulk}}(r):= (r(1-r))^{1+\sigma}, \quad \sigma >0.
	    	\end{equation}
        \end{itemize}
	    	Then, for any fixed $0<x_3\leq 1$, we have 
	    	\begin{equation}\label{q33-BL}
				q_{3,3}^{{\rm st}, L}(\xx,\xx)\lesssim\Ra^2 \left(x_3^{2\beta+4} + \eta^{2\beta}x_3^{4}\right).
			\end{equation}
\end{description}
\end{proposition}

\subsection{Passing from the periodic strip to the periodic half-space}
Let $q_S^{{\rm st}, L}$ satisfy the assumptions of Proposition \ref{Pr1}.
In this section we extend the periodic-strip covariance by zero to the periodic positive half-space. The following lemma shows that this extension satisfies the bounds required for the subsequent half-space estimates.
\begin{lemma}[Extension lemma]\label{extension-lemma}
Suppose $q_S^{{\rm st}, L}$ is a symmetric, positive semidefinite, horizontally stationary covariance kernel on $\DsL\times \DsL$, satisfying the assumptions in Proposition \ref{Pr1}. On $\DhL\times \DhL$ define
\begin{equation*}
    q_S^{{\rm hs}, L}(\yy,\zz):=
    \begin{cases}
        q_S^{{\rm st}, L}(\yy,\zz) &  \yy,\zz\in \DsL\\0 & \mbox{ otherwise}\,, 
    \end{cases}
\end{equation*}
then, $q_S^{{\rm hs}, L}$ is a symmetric, positive semidefinite and horizontally stationary covariance kernel. Moreover we have the following:
\begin{description}
 \item[Case 1]\emph{(Bulk-dominated):} Let $N>2$ and $\gamma>0$. If 
    \[\left|q_S^{{\rm st}, L}(\yy,\zz)\right|\leq (y_3(1-y_3)z_3(1-z_3))^{\gamma}\sum_{\kk'\in\mathbb Z^2}(1+|\yy'-\zz'+L\kk'|)^{-N}\]
    then, for every $\varepsilon>0$ we have
    \[\left|q_S^{{\rm hs}, L}(\yy,\zz)\right|\lesssim m^{\rm hs}(y_3)m^{\rm hs}(z_3)\sum_{\kk'\in\mathbb Z^2}(1+|\yy'-\zz'+L\kk'|)^{-N}\,,\]
    where
    \[m^{\rm hs}(r):=r^{\gamma}(1+r)^{-3-\varepsilon}\]
    and the implicit constant depends only on $\varepsilon$. In particular, choosing $\varepsilon>\gamma-3$,
    $q^{{\rm hs}, L}_S$ satisfies the assumptions of Lemma~\ref{lemma:bulkCov}.
     \item[Case 2]\emph{(Boundary-layer-dominated):} Let $\gamma >0, 0<\eta \leq \frac{1}{2}$, and suppose
     \[q_S^{{\rm st},L}=q_S^{\rm bl,bottom,\eta}+q_S^{\rm bulk,\eta}+q_S^{\rm bl, top,\eta}\]
     satisfies the support and covariance bounds of Proposition \ref{Pr1}. Then $q_S^{{\rm hs},L}$ admits a decomposition 
     \[q_S^{{\rm hs}, L}=q_S^{\rm hs, bl, \eta}+q_S^{\rm hs, bulk, \eta}\]
     into symmetric, positive-semidefinite, horizontally stationary kernels, with the first supported on $y_3,z_3\leq \eta$ and the second supported on $y_3,z_3\geq \eta$, satisfying
      \begin{equation}\label{est1:q_S_hs}
     \left|q_S^{\rm hs,bl,\eta}(\yy,\zz)\right|\lesssim (y_3z_3)^{\gamma}\sum_{\kk'\in \ZZ^2}\left(1+\frac{|\yy'-\zz'+L\kk'|}{\sqrt{y_3 z_3}}\right)^{-N_{\rm bl}^{\rm hs}}
     \end{equation}
     and, for every $\varepsilon>0$,
       \begin{equation}\label{est2:q_S_hs}
       \left|q_S^{\rm hs,bulk,\eta}(\yy,\zz)\right|\lesssim \eta^{2\gamma}m_{\rm bulk}^{\rm hs}(y_3)m_{\rm bulk}^{\rm hs}(z_3)\sum_{\kk'\in \ZZ^2}(1+|\yy'-\zz'+L\kk'|)^{-N_{\rm bulk}^{\rm hs}}\,,
       \end{equation}
       where
       \[m_{\rm bulk}^{\rm hs}(r):=r^{1+\sigma}(1+r)^{-3-\varepsilon}\,\]
       and \[
      N^{\rm hs}_{\rm bl}:=N^{\rm bot}_{\rm bl}>2,\qquad
     N^{\rm hs}_{\rm bulk}:=\min\{N^{\rm bot}_{\rm bl},N_{\rm bulk}\}>2 .
     \]
       The implicit constants in \eqref{est1:q_S_hs} and \eqref{est2:q_S_hs} depend only on $\sigma$ and $\varepsilon$. 
       In particular, choosing $\varepsilon>\sigma-2$, $q^{{\rm hs}, L}_S$ satisfies
    the assumptions of Lemma~\ref{lemma:BL-bulk} with $N_{\rm bl}=N^{\rm hs}_{\rm bl}$ and
    $N_{\rm bulk}=N^{\rm hs}_{\rm bulk}$ and $\gamma = \beta$. 
\end{description}
\end{lemma}

\begin{proof}

Denote $\mathcal{E}q_S$ the extension by zero of a symmetric, horizontally stationary and positive semidefinite covariance function $q_S$ to $\D^{{\rm hs}, L}$. Then clearly $\mathcal{E}q_S$ is also symmetric, horizontal stationary and
\[\int_{\Dh}\int_{\Dh} \overline{f(\yy)}\mathcal{E}q_S(\yy,\zz)f(\zz)\, {\rm d}\yy\, {\rm d}\zz=\int_{\DsL}\int_{\DsL} \overline{f(\yy)}q_S(\yy,\zz)f(\zz)\, {\rm d}\yy\, {\rm d}\zz\geq 0\,.\]

Case 1: We consider $\yy,\zz\in \DsL$ (since outside the extension it vanishes). For $0<r<1$,
\[r^{\gamma}(1-r)^{\gamma} \leq 2^{3+\varepsilon}r^{\gamma}(1+r)^{-3-\varepsilon}= 2^{3+\varepsilon}m^{\rm hs}(r)\,.\]

Consequently,
\[\left|q_S^{{\rm hs}, L}(\yy,\zz)\right|\lesssim m^{\rm hs}(y_3)m^{\rm hs}(z_3) \sum_{\kk'\in \mathbb{Z}^2}(1+|\yy'-\zz'+L\kk'|)^{-N}\,,\]
where the implicit constant is independent of $L$. The same applies outside of $\DsL\times\DsL$ where the extension is zero.

Case 2: Define 
\[q_S^{\rm hs, bl,\eta}:=\mathcal{E}q_S^{\rm bl,bottom,\eta}\,,\quad q_S^{\rm hs, bulk,\eta}:=\mathcal{E}q_S^{\rm bulk,\eta}+\mathcal{E}q_S^{\rm bl,top,\eta} \,.\]

The bottom-layer bound follows directly from the assumption on $q^{\rm bl,bottom,\eta}_S$, since
$N^{\rm hs}_{\rm bl}=N^{\rm top}_{\rm bl}$, and since outside of
$\DsL\times \DsL$ its zero extension vanishes. 
It only remains to estimate the bulk component. We want to prove that 
\[\left|q_S^{\rm hs,bulk, \eta}(\yy,\zz)\right|\lesssim \eta^{2\gamma}m_{\rm bulk}^{\rm hs}(y_3)m_{\rm bulk}^{\rm hs}(z_3)\sum_{\kk'\in \ZZ^2}(1+|\yy'-\zz'+L\kk'|)^{-N_{\rm bulk}^{\rm hs}}\,.\]
Throughout, we use that $t\mapsto(1+t)^{-M}$ is non-increasing in $M$ for $t\geq 0$, and that all
inequalities below hold term by term in the sum over $\kk'\in\ZZ^2$.
We start with the part involving $q_S^{\rm bulk, \eta}$. By assumption on the strip we have 
\[\left|q_S^{\rm bulk, \eta}(\yy,\zz)\right|\leq \eta^{2\gamma}[y_3(1-y_3)]^{1+\sigma}[z_3(1-z_3)]^{1+\sigma}\sum_{\kk'\in \ZZ^2}(1+|\yy'-\zz'+L\kk'|)^{-N_{\rm bulk}}\,.\]
For every $0<r<1$,
\[[r(1-r)]^{1+\sigma}\leq r^{1+\sigma}=(1+r)^{3+\varepsilon}r^{1+\sigma}(1+r)^{-3-\varepsilon}=(1+r)^{3+\varepsilon}m_{\rm bulk}^{\rm hs}(r)\leq 2^{3+\varepsilon}m_{\rm bulk}^{\rm hs}(r)\,, \]
hence we obtain, using also $N_{\rm bulk}^{\rm hs}\leq N_{\rm bulk}$
\[\left|q_S^{\rm bulk, \eta}(\yy,\zz)\right|\leq 2^{6+2\varepsilon}\eta^{2\gamma}m_{\rm bulk}^{\rm hs}(y_3)m_{\rm bulk}^{\rm hs}(z_3)  \sum_{\kk'\in \ZZ^2}(1+|\yy'-\zz'+L\kk'|)^{-N_{\rm bulk}^{\rm hs}}\,.\]
Since outside of $\DsL\times \DsL$ its zero extension vanishes, we proved 
\[\left|\mathcal{E} q_S^{\rm bulk, \eta}(\yy,\zz)\right|\leq  2^{6+2\varepsilon}\eta^{2\gamma}m_{\rm bulk}^{\rm hs}(y_3)m_{\rm bulk}^{\rm hs}(z_3) \sum_{\kk'\in \ZZ^2} (1+|\yy'-\zz'+L\kk'|)^{-N_{\rm bulk}^{\rm hs}}\,.\]
Now we estimate the top layer. The component there vanishes unless $1-\eta\leq y_3,z_3<1$. The assumption on the strip is 
\[\left|q_S^{\rm bl, top, \eta}(\yy,\zz)\right|\leq (1-y_3)^{\gamma}(1-z_3)^{\gamma}\sum_{\kk'\in \ZZ^2}\left(1+\frac{|\yy'-\zz'+L\kk'|}{\sqrt{(1-y_3)(1-z_3)}}\right)^{-N_{\rm bl}^{\rm top}}\,.\]
Since both $1-y_3$ and $1-z_3$ are at most $\eta$, we have 
\[ (1-y_3)^{\beta}(1-z_3)^{\beta}\leq \eta^{2\beta}\,.\]
Further, since $\sqrt{(1-y_3)(1-z_3)}\leq \eta<1$, it holds for every $\kk'\in \ZZ^2$
\[\left(1+\frac{|\yy'-\zz'+L\kk'|}{\sqrt{(1-y_3)(1-z_3)}}\right)^{-N_{\rm bl}^{\rm top}}\leq \left(1+|\yy'-\zz'+L\kk'|\right)^{-N_{\rm bl}^{\rm top}}\leq \left(1+|\yy'-\zz'+L\kk'|\right)^{-N_{\rm bulk}^{\rm hs}}\,, \]
for $N_{\rm bulk}^{\rm hs}\leq N_{\rm bl}^{\rm top}$. We therefore showed that 
\[\left|q_S^{\rm bl, top, \eta}(\yy,\zz)\right|\leq \eta^{2\beta}\sum_{\kk'\in \ZZ^2}\left(1+|\yy'-\zz'+L\kk'|\right)^{-N_{\rm bulk}^{\rm hs}}\,.\]
Now, it is left to reintroduce the weights in the vertical variables. When $\eta\leq \frac 12$, $\frac 12\leq y_3,z_3< 1$, and 
\[m_{\rm{bulk}}^{\rm hs}(r)=\frac{r^{1+\sigma}}{(1+r)^{3+\varepsilon}}\geq \frac{(1/2)^{1+\sigma}}{2^{3+\varepsilon}}=2^{-4-\sigma-\varepsilon}\,,\]
hence 
\[1\leq 2^{8+2\sigma+2\varepsilon}m_{\rm{bulk}}^{\rm hs}(y_3)m_{\rm{bulk}}^{\rm hs}(z_3)\,.\]
In conclusion we showed
\[\left|\mathcal{E}q_S^{\rm bl, top, \eta}(\yy,\zz)\right|\lesssim \eta^{2\beta}m_{\rm{bulk}}^{\rm hs}(y_3)m_{\rm{bulk}}^{\rm hs}(z_3)\sum_{\kk'\in \ZZ^2}\left(1+|\yy'-\zz'+L\kk'|\right)^{-N_{\rm bulk}^{\rm hs}}\,,\]
and outside its support, this inequality is automatically satisfied because the left-hand side is zero.
Finally, we have 
\begin{align*}
    \left|q_S^{\rm hs, bulk, \eta}(\yy,\zz)\right|&=\left|\mathcal{E}q_S^{\rm bl,top, \eta}(\yy,\zz)+\mathcal{E}q_S^{\rm bulk, \eta}(\yy,\zz)\right|\\
    &\leq\left|\mathcal{E}q_S^{\rm bl,top, \eta}(\yy,\zz)\right|+\left|\mathcal{E}q_S^{\rm bulk, \eta}(\yy,\zz)\right|\\
    &\lesssim \eta^{2\beta}m_{\rm{bulk}}^{\rm hs}(y_3)m_{\rm{bulk}}^{\rm hs}(z_3)\sum_{\kk'\in \ZZ^2}\left(1+|\yy'-\zz'+L\kk'|\right)^{-N_{\rm bulk}^{\rm hs}}\,.
\end{align*}
where the implicit constant depends only on $\sigma$ and $\varepsilon$; in particular, it is independent of $L$ and $\eta$.

\end{proof}

\subsection{The unfolding argument}
In the following lemma we identify $\TT_L^2$ with its fundamental cell $[0,L)^2$; translated cells are understood as subsets of $\R^2$. Furthermore we observe the following: fix $\xx\in\DsL$, then
    \[K_{3,3}^{{\rm st}, L}(\xx,\yy)=\sum_{\kk'\in \ZZ^2}K_{3,3}^{\rm st}(\xx;\yy'+L\kk',y_3)\,,\]
    with absolute convergence for almost every $\yy \in\DsL$, with $\xx\neq\yy$.\footnote{Indeed, let $\kk'_0\in\ZZ^2$ be such that $\xx'-\yy'-L\kk'_0\in[-\tfrac L2,\tfrac L2)^2$. For $\kk'\neq \kk'_0$ we have
    $|\xx'-\yy'-L(\kk'-\kk'_0)|\geq L|\kk'-\kk'_0|-|\xx'-\yy'-L\kk'_0|\geq L|\kk'-\kk'_0|-\frac{L}{\sqrt{2}}$. Since $|\kk'-\kk'_0|\geq 1$, we have  $|\xx'-\yy'-L\kk'|=|\xx'-\yy'-L\kk'_0-L(\kk'-\kk'_0)|\geq cL|\kk'-\kk'_0|$ with $c=1-\frac1{\sqrt2}$, so for all but finitely many $\kk'$, the far-field estimate gives
    $|K^{\mathrm{st}}_{3,3}(\xx;\yy'+L\kk',y_3)|\lesssim L^{-3}|\kk'-\kk'_0|^{-3}$, which is summable over $\ZZ^2$; the remaining finitely many terms are finite for $\yy\ne \xx$.
    The single term $\kk'=\kk'_0$ is finite for $\yy\neq \xx$.}
    
\begin{lemma}[Unfolding lemma]\label{unfolding-lemma}
    Let $L\geq 1$. Let $q_S^{{\rm st}, L}$ be the covariance kernel on $\DsL\times \DsL$, identified with its horizontally periodic extension to $\Dis\times \Dis$. Thus 
    \begin{equation}\label{star}
    q_S^{{\rm st}, L}(\yy+L\kk', \zz+L\el')=q_S^{{\rm st}, L}(\yy,\zz), \quad \kk',\el'\in \ZZ^2,\
    \end{equation}
    where horizontal shifts have zero third components.
    Assume that 
    \begin{equation}\label{finite-cov}
    \int_{\Dis}\int_{\Dis}\left|K_{3,3}^{\rm st}(\xx,\yy)\right|\left|q_S^{{\rm st}, L}(\yy,\zz)\right|\left|K_{3,3}^{\rm st}(\xx,\zz)\right|\, {\rm d}\yy\, {\rm d}\zz<\infty\,.
    \end{equation}
    Then the covariance 
    \[q_{3,3}^{{\rm st}, L}(\xx,\xx):=\Ra^2 \int_{\DsL}\int_{\DsL}K_{3,3}^{{\rm st}, L}(\xx,\yy)q_S^{{\rm st}, L}(\yy,\zz) K_{3,3}^{{\rm st}, L}(\xx,\zz)\, {\rm d}\yy\, {\rm d}\zz\]
    admits the unfolded representation
    \begin{equation*}
        q_{3,3}^{{\rm st}, L}(\xx,\xx)=\Ra^2 \int_{\Dis}\int_{\Dis}K_{3,3}^{\rm st}(\xx,\yy)q_S^{{\rm st}, L}(\yy,\zz) K_{3,3}^{\rm st}(\xx,\zz)\, {\rm d}\yy\, {\rm d}\zz\,.
    \end{equation*}
\end{lemma}
\begin{remark}\label{unfolding-remark}
The same statement and proof apply to the half-space,
with $K_{3,3}^{{\rm st},L}$ and
$K_{3,3}^{\rm st}$ replaced by
$K_{3,3}^{{\rm hs},L}$ and
$K_{3,3}^{\rm hs}$, respectively, provided that
$K^{\mathrm{hs},L}_{3,3}(\xx,\yy)=\sum_{\kk'\in\ZZ^2}K^{\rm hs}_{3,3}(\xx;\yy'+L\kk',y_3)$ and that \eqref{finite-cov}
holds with $\Dis$ and $K^{\rm st}_{3,3}$ replaced by $\Dh$ and $K^{\rm hs}_{3,3}$.
Indeed, the unfolding acts only on the horizontal variables. The absolute convergence of the series
follows from the far-field estimate \eqref{eq:far field-new} as in the strip case, and the integrability
condition is established by the proofs of Lemmas~\ref{lemma:bulkCov} and~\ref{lemma:BL-bulk}.
\end{remark}

\begin{remark}
Lemma \ref{unfolding-lemma} uses only the periodicity \eqref{star} of $q_S$; neither symmetry, positive
semidefiniteness nor horizontal stationarity is needed. For horizontally stationary kernels,
$q_S^{{\rm st}, L}(\yy,\zz)=Q(\yy'-\zz',y_3,z_3)$, \eqref{star} is equivalent to the $L$-periodicity of $Q$ in its first argument,
which is the form used in the covariance bounds below. Moreover, in the half-space setting the
integrability assumption \eqref{finite-cov} need not be checked separately: the proofs of
Lemmas~\ref{lemma:bulkCov} and~\ref{lemma:BL-bulk} bound precisely the integral of
$|K_{3,3}|\,|q_S^{{\rm st}, L}|\,|K_{3,3}|$, and hence establish it.
\end{remark}
\begin{proof}
    Denoting 
    \[\DsL+L\kk':=\{(\yy'+L\kk',y_3): (\yy',y_3)\in \DsL\}\,,\]
    where $\DsL$ is identified with $[0,L)^2\times (0,1)$, we first observe that by Tonelli's theorem, using the fact that the integrands are non-negative, we have
    \begin{align*}
    &\int_{\DsL}\int_{\DsL}\left(\sum_{\kk'\in \ZZ^2}|K_{3,3}^{\rm st}(\xx;\yy'+L\kk',y_3)|\right)\left|q_S^{{\rm st}, L}(\yy,\zz)\right|\left(\sum_{\el'\in \ZZ^2}|K_{3,3}^{\rm st}(\xx;\zz'+L\el',z_3)|\right)\,{\rm d}\yy {\rm d}\zz\\
    &=\sum_{\kk',\el'\in \ZZ^2}\int_{\DsL}\int_{\DsL}|K_{3,3}^{\rm st}(\xx;\yy'+L\kk',y_3)|\left|q_S^{{\rm st}, L}(\yy,\zz)\right||K_{3,3}^{\rm st}(\xx;\zz'+L\el',z_3)|\, {\rm d}\yy\, {\rm d}\zz\,.
    \end{align*}
    Changing variables $\yy\rightarrow \yy+L\kk'$ and $\zz\rightarrow \zz+L\el'$ and using \eqref{star} we have
     \begin{align*}
    &\sum_{\kk,\el'\in \ZZ^2}\int_{\DsL+L\el'}\int_{\DsL+L\kk'}|K_{3,3}^{\rm st}(\xx,\yy)|\left|q_S^{{\rm st}, L}(\yy,\zz)\right||K_{3,3}^{\rm st}(\xx,\zz)|\, {\rm d}\yy\, {\rm d}\zz\\
    &=\int_{\Dis}\int_{\Dis}|K_{3,3}^{\rm st}(\xx,\yy)|\left|q_S^{{\rm st}, L}(\yy,\zz)\right||K_{3,3}^{\rm st}(\xx,\zz)|\, {\rm d}\yy\, {\rm d}\zz<\infty\,,
    \end{align*}
    where we used that the translated cells partition $\Dis$ and \eqref{finite-cov}. The preceding calculation establishes absolute integrability. Thus, using the kernel periodisation identity and Fubini's theorem, we obtain
    \begin{align*}
        q_{3,3}^{{\rm st}, L}(\xx,\xx)
        &=\Ra^2 \int_{\DsL}\int_{\DsL}K_{3,3}^{{\rm st}, L}(\xx,\yy)q_S^{{\rm st}, L}(\yy,\zz) K_{3,3}^{{\rm st}, L}(\xx,\zz)\, {\rm d}\yy\, {\rm d}\zz\\
        &=\Ra^2\sum_{\kk',\el'\in \ZZ^2}\int_{\DsL}\int_{\DsL}K_{3,3}^{\rm st}(\xx,\yy+L\kk')q_S^{{\rm st}, L}(\yy,\zz)K_{3,3}^{\rm st}(\xx,\zz+L\el')\, {\rm d}\yy\, {\rm d}\zz\\
        &=\Ra^2\sum_{\kk',\el'\in \ZZ^2}\int_{\DsL+L\el'}\int_{\DsL+L\kk'}K_{3,3}^{\rm st}(\xx,\yy)q_S^{{\rm st}, L}(\yy,\zz)K_{3,3}^{\rm st}(\xx,\zz)\, {\rm d}\yy\, {\rm d}\zz\\
        &=\Ra^2\int_{\Dis}\int_{\Dis}K_{3,3}^{\rm st}(\xx,\yy)q_S^{{\rm st}, L}(\yy,\zz)K_{3,3}^{\rm st}(\xx,\zz)\, {\rm d}\yy\, {\rm d}\zz\,.
    \end{align*}
\end{proof}

\subsection{Covariance bounds on the positive half-space}

In this section we consider the covariance kernel
$q_S^{{\rm hs}, L}$ on
$\DhL\times \DhL$,
obtained by vertical extension-by-zero of
$q_S^{{\rm st}, L}$ in Lemma \ref{extension-lemma}.
We identify $q_S^{{\rm hs}, L}$ with its horizontally
$L$-periodic extension to
$\Dh\times \Dh$, using the same notation for both.
In particular,
\[
    q_S^{{\rm hs}, L}(\yy,\zz)=0
    \qquad\text{if }y_3\geq 1\text{ or }z_3\ge1.
\]

By the half-space version of the unfolding lemma
(see Remark \ref{unfolding-remark}), the periodic half-space velocity covariance
admits the representation
\begin{equation}\label{hs-unfolded-representation}
    q_{3,3}^{{\rm hs}, L}(\xx,\xx)=
    \Ra^2\int_{\Dh}\int_{\Dh}
    K_{3,3}^{\rm hs}(\xx,\yy)\,q_S^{{\rm hs}, L}(\yy,\zz)\,
    K_{3,3}^{\rm hs}(\xx,\zz)
    \,{\rm d}\yy\,{\rm d}\zz\,.
\end{equation}
Here $K_{3,3}^{\rm hs}$ is the nonperiodic half-space
Stokes kernel. The estimates below establish an upper bound for the right-hand side of \eqref{hs-unfolded-representation}.

To simplify notation, throughout this subsection we write
\[K_{3,3}:=K_{3,3}^{\rm hs},
    \qquad
    q_S:=q_S^{{\rm hs}, L}.\]
All subsequent integrals are taken over the
nonperiodic half-space, while the covariance remains
horizontally $L$-periodic.

By horizontal stationarity and horizontal translation
invariance of the Stokes problem, the diagonal velocity
covariance depends only on $x_3$. We set
\[
    \widetilde q_{3,3}(x_3)
    :=
    q_{3,3}^{{\rm hs}, L}
    \left((\bm{0}',x_3),(\bm{0}',x_3)\right).
\]
The dependence of $\widetilde q_{3,3}$ on $L$ is suppressed
in the notation.

\begin{lemma}[Bounds on $\widetilde{q}_{3,3}$, case 1]\label{lemma:bulkCov}
Let $L\geq 1$. Suppose that the covariance kernel satisfies
\begin{equation}\label{bulkCov-scaling-deacy}
|q_S(\yy,\zz)|\lesssim m(y_3)m(z_3)\sum_{\kk'\in \mathbb{Z}^2}(1+|\yy'-\zz'+L\kk'|)^{-N}, \quad N>2\,,
\end{equation}
where 
$$m(r)=r^{\gamma}(1+r)^{-3-\varepsilon},\qquad 0<\gamma<3+\varepsilon, \quad \varepsilon>0\,.$$
Then, for any fixed $0<x_3\leq 1$, we have the upper bound
\begin{equation*}
	\widetilde{q}_{3,3}(x_3)\lesssim \Ra^2 x_3^{4}\,.
\end{equation*}
\end{lemma}
\begin{proof}
Inserting \eqref{bulkCov-scaling-deacy} into \eqref{hs-unfolded-representation} and using Cauchy-Schwarz as well as Young's inequality, we obtain 
\begin{align*}
	\widetilde{q}_{3,3}(x_3)&\lesssim \Ra^2 \int_\HS \int_\HS |K_{3,3}((\bzero',x_3),\yy)|m(y_3)m(z_3)\sum_{\kk'\in\ZZ^2}(1+|\yy'-\zz'+L\kk'|)^{-N}\\
    &\qquad\qquad\qquad\quad|K_{3,3}((\bzero',x_3)),\zz)|\, {\rm d}\yy\,{\rm d}\zz\\
	       &\leq \Ra^2 \int_0^{\infty}\int_0^{\infty}\|f_{x_3}(\cdot, y_3)\|_{L^1(\R^2)}\|h^L\|_{L^{\infty}(\R^2)}\|f_{x_3}(\cdot, z_3)\|_{L^1(\R^2)}\, {\rm d}y_3\, {\rm d}z_3\,,
\end{align*}
where
\[f_{x_3}(\cdot,y) := m(y)K_{3,3}((\bm{0}',x_3),(\cdot,y)), \quad h^L(\cdot) := \sum_{\kk'\in\mathbb{Z}^2}(1+ |\cdot+L\kk'|)^{-N}.\]
Since by assumption it holds $N>2$, we have\footnote{Indeed
\[
\|h^L\|_{L^\infty(\mathbb R^2)}
=h^L(\bm{0}')
=\sum_{\kk'\in\mathbb Z^2}(1+L|\kk'|)^{-N}
\lesssim_N 1+L^{-2},
\quad N>2.\]} 
\[\|h^L\|_{L^{\infty}}\lesssim_N 1+L^{-2}\,.\]

we obtain 
\begin{equation}\label{q33BoundmK}
\widetilde{q}_{3,3}(x_3)\lesssim \Ra^2\left(1+\frac{1}{L^2}\right)\left(\int_0^{\infty}m(y_3)\int_{\R^2}|K_{3,3}((\bzero',x_3), (\yy',y_3))|\,{\rm d}\yy'\,{\rm d}y_3\right)^2\,.
\end{equation}
In order to bound the right-hand side, we split
$(0,\infty)=(0,\frac{x_3}{9})\cup [\frac{x_3}{9}, 9x_3]\cup (9x_3, \infty)$.\\
\\

\underline{Region $y_3\in (0,\frac{x_3}{9})$.}
In this region we have $\min\{x_3,y_3\}=y_3$, $x_3-y_3\sim x_3$ and
$$|\xx-\yy|\geq |x_3-y_3|=x_3-y_3\geq 8y_3=8\min\{x_3,y_3\}.$$
Therefore we use the far field estimate \eqref{eq:far field-new} and obtain
$$|K_{3,3}((\bzero', x_3),\yy)|\lesssim \frac{\min\{x_3,y_3\}^2}{(|\yy'|^2+x_3^2)^{\frac 32}}\sim \frac{y_3^2}{(|\yy'|^2+x_3^2)^{\frac 32}}\,.$$
Integrating over $\R^2$ we obtain
\begin{equation*}
\int_{\R^2}|K_{3,3}((\bzero',x_3), (\yy',y_3))|\,{\rm d}\yy'=y_3^2\int_{\R^2}\frac{1}{(|\yy'|^2+x_3^2)^{\frac 32}}{\rm d}\yy'\sim \frac{y_3^2}{x_3}\,,
\end{equation*}
hence, for $\gamma>-3$
$$\int_{0}^{\frac{x_3}{9}}m(y_3)\frac{y_3^2}{x_3}\,{\rm d}y_3 \sim \frac{1}{x_3}\int_{0}^{\frac{x+3}{9}}y_3^{\gamma+2}\, {\rm d}y_3=x_3^{\gamma+2}\,.$$

\underline{Region $y_3\in (\frac{x_3}{9}, 9x_3)$.}
In this region it holds $x_3\sim y_3\sim\min\{x_3,y_3\}$, and 
\begin{equation}\label{BoundKL1}
\int_{\R^2}|K_{3,3}((\bzero',x_3),(\yy',y_3))|\,{\rm d}\yy'\lesssim x_3\,.
\end{equation}
Indeed, we split the integral into  
$$\left(\int_{|\xx-\yy|\leq 8\min\{x_3,y_3\}}+\int_{|\xx-\yy|\geq 8\min\{x_3,y_3\}}\right)|K_{3,3}((\bzero',x_3),(\yy',y_3))|\,{\rm d}\yy'\,.$$ 
In the near field, the condition $\sqrt{|\yy'|^2+(x_3-y_3)^2}\leq 8\min\{x_3,y_3\}$ implies 
$$|\yy'|\leq \sqrt{64\min\{x_3,y_3\}^2-(x_3-y_3)^2}\leq 8 \min\{x_3,y_3\}\,.$$
 Hence, using the near-field kernel estimate \eqref{eq:near-field} we have
\begin{align*}
&\int_{|\xx-\yy|\leq 8\min\{x_3,y_3\}} |K_{3,3}((\bzero',x_3),(\yy',y_3))|\,{\rm d}\yy'\\
&\lesssim \int_{|\yy'|\leq 8 \min\{x_3,y_3\}}\frac{1}{\sqrt{|\yy'|^2+(x_3-y_3)^2}}\,{\rm d}\yy'\\
&=2\pi \int_0^{8\min\{x_3,y_3\}}\frac{\rho}{\sqrt{\rho^2+(x_3-y_3)^2}} \, {\rm d}\rho\\
&\leq 2\pi  \int_0^{8\min\{x_3,y_3\}} 1\, {\rm d}\rho\lesssim x_3\,.
\end{align*}
In the far field, we use the far-field kernel estimate \eqref{eq:far field-new} and, passing to polar coordinates, we obtain
\begin{align*}
    &\int_{|\xx-\yy|\geq 8\min\{x_3,y_3\}} |K_{3,3}((\bzero',x_3),(\yy',y_3))|\,{\rm d}\yy'\\
    &\lesssim\min\{x_3,y_3\}^2\int_{\sqrt{64\min\{x_3,y_3\}^2-(x_3-y_3)^2}}^{\infty}\frac{\rho}{(\rho^2+(x_3-y_3)^2)^{\frac 32}}\, {\rm d}\rho\,.
\end{align*}
Notice that the lower limit is well defined because in this region $|x_3-y_3|<8\min\{x_3,y_3\}$.
Setting $s=\rho^2+(x_3-y_3)^2, \; {\rm d}s = 2\rho\, {\rm d}\rho$, we get
\begin{align*}
 &\min\{x_3,y_3\}^2\int_{\sqrt{64\min\{x_3,y_3\}^2-(x_3-y_3)^2}}^{\infty}\frac{\rho}{(\rho^2+(x_3-y_3)^2)^{\frac 32}}\, {\rm d} \rho\\
&=\min\{x_3,y_3\}^2\int_{64\min\{x_3,y_3\}^2}^{\infty}\frac{1}{2s^{\frac 32}}{\rm d}s=\frac{\min\{x_3,y_3\}}{8} \lesssim x_3\,.
 \end{align*}
 This proves \eqref{BoundKL1}. Therefore 
 \begin{equation*}
 \int_{\frac{x_3}{9}}^{9x_3}m(y_3)\|K_{3,3}((\bzero',x_3),(\cdot, y_3))\|_{L^1(\R^2)}\, {\rm d}y_3\lesssim x_3m(x_3) \int_{\frac{x_3}{9}}^{9x_3}\, {\rm d}y_3\lesssim  x_3^{\gamma+2}.
 \end{equation*}

\underline{Region $y_3\in (9x_3,\infty)$.}
In this region we have $\min\{x_3,y_3\}=x_3$, $y_3-x_3\sim y_3$ and 
$$|\xx-\yy|\geq |x_3-y_3|=y_3-x_3\geq 8x_3=8\min\{x_3,y_3\}.$$
Therefore we use the far-field kernel estimate \eqref{eq:far field-new} and obtain
$$|K_{3,3}((\bzero', x_3),(\yy',y_3))|\lesssim \frac{x_3^2}{(|\yy'|^2+y_3^2)^{\frac 32}}\,.$$
Integrating in $\R^2$ we obtain
\begin{equation*}
\int_{\R^2}|K_{3,3}((\bzero',x_3),(\yy',y_3))|\,{\rm d}\yy'=x_3^2\int_{\R^2}\frac{1}{(|\yy'|^2+y_3^2)^{\frac 32}}\,{\rm d}\yy'\sim \frac{x_3^2}{y_3}\,,
\end{equation*}
hence it holds
\begin{align*}
    \int_{9x_3}^{\infty}m(y_3)\|K_{3,3}((\bzero',x_3),(\yy',y_3))\|_{L^1(\R^2)}\,{\rm d}y_3&\lesssim \int_{9x_3}^{\infty}m(y_3)\frac{x_3^2}{y_3}\,{\rm d}y_3 \\
    &\lesssim x_3^2\int_{9x_3}^{\infty}y_3^{\gamma-1}(1+y_3)^{-3-\varepsilon}\,{\rm d}y_3\,.
\end{align*}
Now, if $0<x_3<\frac 19$, we split the integral on the right-hand side
\begin{align*}
&\int_{9x_3}^{\infty}y_3^{\gamma}(1+y_3)^{-3-\varepsilon}\frac{1}{y_3^2}\,{\rm d}y_3\\
&=\left(\int_{9x_3}^{1}+\int_{1}^{\infty}\right)y_3^{\gamma-1}(1+y_3)^{-3-\varepsilon}\,{\rm d}y_3\\
&\leq \int_{9x_3}^{1} y_3^{\gamma-1}\, {\rm d}y_3+\int_{1}^{\infty} y_3^{\gamma-4-\varepsilon}\, {\rm d}y_3 
\end{align*}
The first integral is 
\begin{equation}\label{far-field-gamma}
\int_{9x_3}^{1} y_3^{\gamma-1}\, {\rm d}y_3=\frac{1-(9x_3)^{\gamma}}{\gamma} \quad \mbox{ for }\gamma\neq 0\,.
\end{equation}
The second integral is finite thanks to the assumption $\gamma<3+\varepsilon$.
Consequently, recalling that we are only interested in $\gamma>0$ we have 
\[x_3^2\int_{9x_3}^{\infty}y_3^{\gamma-1}(1+y_3)^{-3-\varepsilon}\,{\rm d}y_3\lesssim x_3^2\,.\]

If $\frac 19\leq x_3\leq 1$, we estimate directly
\[\int_{9x_3}^{\infty} y_3^{\gamma-1}(1+y_3)^{-3-\varepsilon}\, {\rm d}y_3\leq x_3^2\int_1^{\infty} y_3^{\gamma-4-\varepsilon}\, {\rm d}y_3\lesssim x_3^2\,.\]
Combining all the estimates in these three regions we obtain, for $0< x_3\leq 1$ 
\[\int_{9x_3}^{\infty}m(y_3)\|K_{3,3}((\bzero',x_3),(\cdot,y_3))\|_{L^1(\R^2)}\,{\rm d}y_3\lesssim x_3^2 \qquad \mbox{ for } 0<\gamma<3+\varepsilon\,.\]
It is left to finally combine the estimates we obtained in the regions 
$(0,\frac{x_3}{9}), [\frac{x_3}{9}, 9x_3], (9x_3, \infty)$. We obtained
\begin{equation}\label{BoundKmL1}
  \int_0^{\infty}m(y_3)\int_{\R^2}|K_{3,3}((\bzero',x_3), (\yy',y_3))|\,{\rm d}\yy'\,{\rm d}y_3\lesssim x_3^{\beta+2}+x_3^2\,,   
\end{equation}
from which we conclude
\begin{equation*}
\widetilde{q}_{3,3}(x_3)\lesssim \Ra^2\left(1+\frac{1}{L^2}\right)\left(x_3^{\beta+2}+x_3^2\right)^2\lesssim \Ra^2\left(x_3^{\gamma+2}+x_3^2\right)^2\,,
\end{equation*}
where in the last inequality we used that $L\geq 1$ by assumption.
\end{proof}

\begin{lemma}[Bounds on $\widetilde{q}_{3,3}$, case 2]\label{lemma:BL-bulk}
Let $L\geq 1$ and $0<\eta<1$. Suppose that the covariance kernel decomposes into a boundary-layer component and a bulk component:
$$q_S=q_S^{\mathrm{bl},\eta}+q_S^{\mathrm{bulk},\eta}$$
and each component is horizontally stationary, thus we may write 
$$q_S^\ast(\yy,\zz)=q_S^\ast(\yy'-\zz',y_3,z_3),\qquad\ast\in\{\mathrm{bl},\mathrm{bulk}\}.$$

Assume that the boundary-layer covariance $q_S^{\mathrm{bl},\eta}$ is localized near the wall:
\begin{equation*}
	q_S^{\mathrm{bl},\eta}(\yy'-\zz', y_3, z_3)=0\qquad \text{unless } y_3,z_3\leq \eta.
\end{equation*}

Moreover, we assume for $\gamma >0$
\begin{equation}\label{scaling-decay:bl}
\left|q_S^{\rm{bl},\eta}(\yy'-\zz',y_3,z_3)\right|\lesssim (y_3z_3)^\gamma\sum_{\kk'\in \ZZ^2}\left(1+\frac{|\yy'-\zz'+L\kk'|}{\sqrt{y_3z_3}}\right)^{-N_{\rm{bl}}},\qquad N_{\rm{bl}}>2.
\end{equation}

We assume that the bulk covariance $q_S^{\mathrm{bulk},\eta}$ is supported away from the wall:
\begin{equation*}
q_S^{\mathrm{bulk},\eta}(\yy'-\zz', y_3, z_3)=0\qquad\text{unless } y_3,z_3 > \eta.
\end{equation*}

Furthermore, we assume
\begin{equation}\label{scaling-decay:bulk}
\left|q_S^{\rm{bulk},\eta}(\yy'-\zz',y_3,z_3)\right|\lesssim\eta^{2\gamma}\,m_{\rm{bulk}}(y_3)m_{\rm{bulk}}(z_3)\sum_{\kk'\in \ZZ^2}\left(1+|\yy'-\zz'+L\kk'|\right)^{-N_{\rm{bulk}}},\, N_{\rm{bulk}}>2,
\end{equation}
where the bulk profile satisfies
$$m_{\mathrm{bulk}}(r):=r^{1+\sigma}(1+r)^{-3-\varepsilon},\qquad 0<\sigma<2+\varepsilon, \quad \varepsilon>0.$$
Then, for any fixed $0<x_3\leq 1$, we have the upper bound
\begin{equation}\label{q33-Malkus}
	\widetilde{q}_{3,3}(x_3)\lesssim \Ra^2 (x_3^{2\gamma+4}+\eta^{2\gamma}x_3^4)\,.
\end{equation}
The implicit constant depends only on $\gamma$, $N_{\rm bl}$, $N_{\rm bulk}$, $\sigma$ and $\varepsilon$; in particular, it is uniform in $L\geq 1$ and $0<\eta<1$.
\end{lemma}

\begin{proof}

\textbf{Boundary-layer covariance:}
Using the assumption \eqref{scaling-decay:bl} in the representation
	$$\widetilde{q}_{3,3}^{\rm{bl}, \eta}(x_3) :=\Ra^2\int_\HS\int_\HS K_{3,3}((\bzero',x_3), \yy)q_S^{\rm{bl},\eta}(\yy,\zz)K_{3,3}((\bm{0}',x_3), \zz)\, {\rm d}\yy\, {\rm d}\zz $$
	we obtain 
	\begin{align*}	
		&\left|\widetilde{q}_{3,3}^{\rm{bl},\eta}(x_3)\right|\lesssim\\ &\Ra^2\int_{0}^{\eta}\int_{0}^{\eta}\int_{\R^2}\int_{\R^2} |K_{3,3}((\bzero',x_3), (\yy',y_3))|(y_3z_3)^\beta\sum_{\kk'\in \ZZ^2}\left(1+\frac{|\yy'-\zz'+L\kk'|}{\sqrt{y_3z_3}}\right)^{-N_{\rm{bl}}}\\
        &\qquad\qquad\qquad\qquad\quad |K_{3,3}((\bzero',x_3), (\zz',z_3))|\,{\rm d}\yy'\, {\rm d}\zz'\, {\rm d}y_3\, {\rm d}z_3\\
        &\leq \Ra^2 \int_0^{\eta} \int_0^{\eta} y_3^\gamma\|f_{x_3}(\cdot,y_3)\|_{L^1(\mathbb{R}^2)}\|h^L_{y_3,z_3}\|_{L^{\infty}(\mathbb{R}^2)}z_3^\gamma\|f_{x_3}(\cdot, z_3)\|_{L^1(\mathbb{R}^2)}{\rm d}y_3{\rm d}z_3
    \end{align*}
    where
    \[f_{x_3}(\cdot,y) := K_{3,3}((\bm{0}',x_3),(\cdot,y)), \quad h^L_{y_3,z_3}(\xib') := \sum_{\kk'\in \ZZ^2}\left(1+ \frac{|\xib'+L\kk'|}{\sqrt{y_3 z_3}}\right)^{-N_{\rm{bl}}}.\]
    Observe that due to $N_{\rm bl} >2$ and $\sqrt{y_3z_3}\leq \eta<1$ it holds\footnote{Indeed
\[
\|h^L_{y_3,z_3}\|_{L^\infty(\mathbb R^2)}
=h^L_{y_3,z_3}(\bm{0}')
=\sum_{\kk'\in\mathbb Z^2}\left(1+L\frac{|\kk'|}{\sqrt{y_3z_3}}\right)^{-N}
\lesssim_N 1+\frac{y_3z_3}{L^2} \lesssim_N 1,
\quad N>2.\]}
    \[\|h_{y_3,z_3}^L\|_{L^{\infty}(\R^2)}\lesssim 1\,.\]
    Hence we obtain
	\begin{align*}
		\left|\widetilde{q}_{3,3}^{\rm{bl}, \eta}(x_3)\right|
		\lesssim \Ra^2\left(\int_0^{\eta}y_3^{\gamma}\|f_{x_3}(\cdot,y_3)\|_{L^1(\R^2)}\, {\rm d}y_3\right)^2=:\Ra^2 I_{\eta}^2(x_3)\,.
	\end{align*}
	We decompose\footnote{ The decomposition is chosen according to the fact that the near field condition $|\xx-\yy|\leq M$ implies  $|x_3-y_3|\leq M$. If $\min\{x_3,y_3\}=x_3$ then $y_3\leq 9x_3$, while if $\min\{x_3,y_3\}=y_3$, then $x_3\leq 9y_3$. } 
	the $y_3$-domain $(0,\eta)\subset(0, {\frac{x_3}{9}})\cup ({\frac{x_3}{9}},{9x_3})\cup ({9x_3}, {\eta})$. Denote $\xx=(\bm{0}',x_3)$, $\yy=(\yy',y_3)$, $M:=8\min\{x_3,y_3\}$.
    
    \underline{Region $y_3\in \left(0,\frac{x_3}{9}\right)$.} 
	Here it holds $\min\{y_3,x_3\}=y_3$ and $|\xx-\yy|\geq x_3-y_3\geq 8y_3=M$
	\begin{align*}
		\|f_{x_3}(\cdot,y_3)\|_{L^1(\R^2)}\lesssim \int_{\R^2}\frac{y_3^2}{|\xx-\yy|^3}\,{\rm d}\yy' &\sim\int_{\R^2} \frac{y_3^2}{(|\yy'|^2+(x_3-y_3)^2)^{\frac 32}}\,{\rm d}\yy'\\
		&\sim y_3^2\int_0^{\infty}\frac{\rho}{(\rho^2+(x_3-y_3)^2)^{\frac 32}}\, {\rm d}\rho\\ &\sim\frac{y_3^2}{(x_3-y_3)}\,.
	\end{align*}
	Since $\frac{8}{9}x_3\leq x_3-y_3\leq x_3$, then $x_3-y_3\sim x_3 $ and
	\begin{equation*}
		\|f_{x_3}(\cdot,y_3)\|_{L^2(\R^2)}\lesssim  \frac{y_3^2}{x_3}\,.
	\end{equation*}
    Consequently
\[
\int_0^{x_3/9} y_3^{\gamma}\,\|f_{x_3}(\cdot,y_3)\|_{L^1(\mathbb{R}^2)}\,\mathrm{d}y_3
\lesssim \frac{1}{x_3}\int_0^{x_3/9} y_3^{\gamma+2}\,\mathrm{d}y_3 \lesssim x_3^{\gamma+2}.
\]
	\underline{Region $y_3 \in \left(\frac{x_3}{9}, 9 x_3\right)$.}  Here $x_3\sim y_3\sim\min\{x_3,y_3\}$. We have $|x_3-y_3|\leq M$, so that $\sqrt{M^2-|x_3-y_3|^2}$ is well-defined and $\{|\xx-\yy|\leq M\}=\{|\yy'|\leq \sqrt{M^2-|x_3-y_3|^2}\}$ \footnote{Notice that in this region, $|x_3-y_3|\leq M$. Indeed if $x_3\leq y_3$ then $y_3\leq 9 x_3$, so $|x_3-y_3|\leq 8\min\{x_3,y_3\}$.
	 If instead $y_3\leq x_3$ then $x_3\leq 9y_3$, so $|x_3-y_3|\leq 8\min\{x_3,y_3\}$.}. We write 
	\begin{equation*}
	\|f_{x_3}(\cdot,y_3)\|_{L^1(\R^2)}\lesssim \int_{|\xx-\yy|\leq M}\frac{1}{|\xx-\yy|}\, {\rm d}\yy'+\int_{|\xx-\yy|\geq M}\frac{M^2}{|\xx-\yy|^3}\, {\rm d}\yy'\,,
	\end{equation*}
	where we used the kernel estimates \eqref{eq:near-field} and \eqref{eq:far field-new}.
	The first integral is
	\begin{align*}\label{argument-nearfield}
		&\int_{|\xx-\yy|\leq M}\frac{1}{|\xx-\yy|}\, {\rm d}\yy'\\
		&\sim \int_0^{\sqrt{M^2-|x_3-y_3|^2}} \frac{\rho }{(\rho^2+(x_3-y_3)^2)^{\frac 12}}\,{\rm d}\rho \\&\sim (\sqrt{M^2-|x_3-y_3|^2+|x_3-y_3|^2}-|x_3-y_3|)\\
		&=(M-|x_3-y_3|)\leq M
	\end{align*}
	The second integral instead is
		\begin{align*}
		&\int_{|\xx-\yy|\geq M}\frac{M^2}{|\xx-\yy|^3}\, {\rm d}\yy'\\
        &\lesssim M^2\int_{\sqrt{M^2-|x_3-y_3|^2}}^{\infty}\frac{\rho}{(\rho^2+|x_3-y_3|^2)^{\frac 32}} {\rm d}\rho\sim \frac{M^2}{\sqrt{M^2-|x_3-y_3|^2+|x_3-y_3|^2}}=M\,.
		\end{align*}
        Since $M=8\min\{x_3,y_3\}\leq 8x_3$, we conclude
        $\|f_{x_3}(\cdot,y_3)\|_{L^1(\mathbb{R}^2)}\lesssim x_3$ and therefore
\[
\int_{x_3/9}^{9x_3} y_3^{\gamma}\,\|f_{x_3}(\cdot,y_3)\|_{L^1(\mathbb{R}^2)}\,\mathrm{d}y_3
\lesssim x_3^{\gamma+2}.
\]
	\underline{Region $y_3 \in \left(9x_3, \eta\right)$.} In this region it holds
    $$|\xx-\yy|\geq |x_3-y_3|=y_3-x_3\geq 8x_3=8\min\{x_3,y_3\}.$$
    Hence here we compute
    \begin{align*}
		\|f_{x_3}(\cdot,y_3)\|_{L^1(\R^2)}
		\lesssim \int_{\R^2}\frac{M^2}{|\xx-\yy|^3}\,{\rm d}\yy'
		&\sim \int_{\R^2}\frac{x_3^2}{(|\yy'|^2+(x_3-y_3)^2)^{\frac 32}}\, {\rm d}\yy'\\ &\sim x_3^2\int_0^{\infty}\frac{\rho}{(\rho^2+(x_3-y_3)^2)^{\frac 32}}\, {\rm d}\rho\\ &\sim \frac{x_3^2}{y_3-x_3}\,.
	\end{align*}
	Since here $\frac{8}{9}y_3\leq y_3-x_3\leq y_3$, then $y_3-x_3\sim y_3$ and 
	\begin{equation*}
		\|f_{x_3}(\cdot,y_3)\|_{L^1(\R^2)}\lesssim \frac{x_3^2}{y_3}\,,
	\end{equation*}
    implying 
    \[\int_{9x_3}^{\eta} y_3^{\gamma}\,\|f_{x_3}(\cdot,y_3)\|_{L^1(\mathbb{R}^2)}\,\mathrm{d}y_3
     \lesssim \frac{x_3^2\,\eta^\gamma}{\gamma}\]

Collecting all estimates, we have
\begin{align*}
I_\eta(x_3)
&=\int_0^\eta y_3^{\beta}\,\|f_{x_3}(\cdot,y_3)\|_{L^1(\mathbb{R}^2)}\,\mathrm{d}y_3\\
&\lesssim x_3^{\gamma+2}+x_3^{\gamma+2}+\frac{x_3^2\,\eta^\gamma}{\gamma},
\end{align*}
where each integral is understood over its interval intersected with $(0,\eta)$, and $\gamma>0$
ensures the bound on the last integral. Hence, we proved
\begin{equation}\label{BLest33}
\widetilde q^{\mathrm{bl},\eta}_{3,3}(x_3)\lesssim \mathrm{Ra}^2 I_\eta^2(x_3)
\lesssim \mathrm{Ra}^2\big(x_3^{2\gamma+4}+x_3^4\eta^{2\gamma}\big).
\end{equation}

	\textbf{Bulk covariance:} Using the assumption \eqref{scaling-decay:bulk} in the representation
	$$\widetilde{q}_{3,3}^{\rm{bulk}, \eta}(x_3) :=\Ra^2\int_\HS\int_\HS K_{3,3}((\bzero',x_3), \yy)q_S^{\rm{bulk},\eta}(\yy,\zz)K_{3,3}((\bzero',x_3), \zz)\, {\rm d}\yy\, {\rm d}\zz $$
	we obtain with an argumentation as before 
	\begin{align*}	
		&\left|\widetilde{q}_{3,3}^{\rm{bulk}, \eta}(x_3)\right|\lesssim\\
        &\Ra^2\eta^{2\beta}\int_{\eta}^{\infty}\int_{\eta}^{\infty}\int_{\R^2}\int_{\R^2} |K_{3,3}((\bzero',x_3), (\yy',y_3)))| m_{\rm{bulk}}(y_3)\sum_{\kk'\in \ZZ^2} (1+|\yy'-\zz'+L\kk'|)^{-N_{\rm bulk}} \\
        &\qquad \qquad \qquad \qquad\qquad \quad m_{\rm{bulk}}(z_3)|K_{3,3}((\bzero',x_3), (\zz',z_3))|\,{\rm d}\yy'\, {\rm d}\zz'\, {\rm d}y_3\, {\rm d}z_3\\
        &\leq \Ra^2\eta^{2\beta}\int_\eta^\infty\,\int_\eta^\infty m_{\rm{bulk}}(y_3)\,\|f_{x_3}(\cdot,y_3)\|_{L^1(\R^2)}\,\|h^L\|_{L^\infty(\mathbb{R}^2)}\,m_{\mathrm{bulk}}(z_3)\,\|f_{x_3}(\cdot,z_3)\|_{L^1(\R^2)}\,\mathrm{d}y_3\,\mathrm{d}z_3,
	\end{align*}
    where $f_{x_3}$ is as before and
    \[
    h^L(\xi'):=\sum_{\kk'\in\ZZ^2}\big(1+|\xi'+L\kk'|\big)^{-N_{\rm{bulk}}} .
    \]
    By the same argument as for $h^L_{y_3,z_3}$ (now with $\sqrt{y_3z_3}$ replaced by $1$), and since
    $N_{\rm{bulk}}>2$, it holds
    \begin{align*}
        \left|\widetilde{q}_{3,3}^{\rm{bulk}, \eta}(x_3)\right|&\lesssim\Ra^2\eta^{2\beta} \left(\int_{\eta}^{\infty}m_{\rm{bulk}}(y_3)\int_{\R^2}|K_{3,3}((\bm{0}',x_3),(\yy',y_3))|\,{\rm d}\yy'\, {\rm d}y_3\right)^2 \notag\\
        &=: \Ra \eta^{2\gamma}I_{\eta}(x_3)^2.
    \end{align*}
    Since the integrand is nonnegative,
    \begin{equation*}
    I_\eta(x_3)\le \int_0^\infty m_{\rm{bulk}}(y_3)\int_{\R^2}|K_{3,3}((\bm{0}',x_3),(\yy',y_3))|\,\mathrm{d}\yy'\,\mathrm{d}y_3 ,
     \end{equation*}
    which is the integral on the right-hand side of \eqref{q33BoundmK} with $m$ replaced by $m_{\rm{bulk}}$, that is, with $\gamma=\sigma+1$. Since $0<\sigma<2+\varepsilon$, we have $0<\gamma< 3+\varepsilon$, so the three-region estimate in the proof of Lemma \ref{lemma:bulkCov} applies, and yields, for $0< x_3\leq 1$\footnote{see estimate \eqref{BoundKmL1}}, 
    \[I_{\eta}(x_3)\lesssim x_3^{\sigma+3}+x_3^2\lesssim x_3^2\,,\]
    and, therefore
    \begin{equation}\label{BulkEst33}\left|\widetilde{q}_{3,3}^{\rm{bulk}, \eta}(x_3)\right|\lesssim \Ra^2 \eta^{2\gamma}x_3^4\,.
    \end{equation}

	Putting together \eqref{BLest33} and \eqref{BulkEst33}, we obtain \eqref{q33-Malkus}.
\end{proof}

\begin{remark}[The role of $\sigma$]\label{rem:sigma}
In the bulk part of the proof, $I_\eta$ is bounded by the three-region estimate of
Lemma \ref{lemma:bulkCov} with $\gamma=\sigma+1$, and $\sigma$ only enters through the far-field region
$y_3>9x_3$, more precisely through the integral \eqref{far-field-gamma}. For $\gamma>0$, i.e.\ $\sigma>-1$,
this integral is bounded by $1/\gamma$, uniformly in $x_3$. For a bulk covariance of order one,
$\gamma=0$, it becomes logarithmic; since the bulk covariance is supported in $y_3>\eta$, the lower
limit of integration is then $\max\{9x_3,\eta\}$, which produces a factor $\log(1/\eta)$ in $I_\eta$, hence
$\log^2(1/\eta)$ in \eqref{BulkEst33}, and logarithmic corrections in the bound on the Nusselt number.
\end{remark}

\subsection{Passing from the positive half space to the strip: corrector estimates and proof of Proposition \ref{Pr1}}\label{sec:corrector}

 Let $L>0$ and recall that we denote $\D^{{\rm st}, L}=\TT^2\times (0,1)$, $\D^{{\rm hs}, L}=\TT^2\times (0,\infty)$, and the horizontal average
\[\langle g\rangle' (x_3):=\frac{1}{L^2}\int_{\TT^2_L}g(\xx',x_3)\, {\rm d}\xx'\]
for any $g$ on $\D^{{\rm st}, L}$. For $\xx\in \D^{{\rm st}, L}$ write 
\begin{equation}\label{K-dec}
    K_{3,3}^{{\rm st}, L}(\xx,\yy)=K_{3,3}^{{\rm hs}, L}(\xx,\yy)-H^L(\xx,\yy)\,,\quad \yy\in \D^{{\rm st}, L}\,,
\end{equation}
where for each fixed $\xx,\yy\mapsto H^L(\xx,\yy)$ is the unique horizontally $L$-periodic solution of 
\begin{equation}\label{HL}
\begin{cases}
\Delta^2_{\yy} H^L(\xx,\cdot)=0 & \mbox{ in }\D^{{\rm st}, L}\\
H^L=\partial_{y_3} H^L=0 & \mbox{ on }\{y_3=0\}\\
H^L(\xx,\cdot)=K_{3,3}^{{\rm hs},L}(\xx,\cdot)\,, \quad \partial_{y_3}H^L(\xx,\cdot)=\partial_{y_3}K_{3,3}^{{\rm hs},L}(\xx,\cdot) & \mbox{ on }\{y_3=1\}.
\end{cases}
\end{equation}
Recall that, for kernels $\Phi,\Psi$ on $\D^{{\rm st},L}$, we set
\[B(\Phi,\Psi):=\Ra^2\int_{\D^{{\rm st},L}}\int_{\D^{{\rm st},L}}\Phi(\yy)\,q^{{\rm st},L}_S(\yy,\zz)\,\Psi(\zz)\,{\rm d}\yy\,{\rm d}\zz.\]
Since $q^{{\rm st},L}_S$ is positive semidefinite, $\|\Phi\|_B:=B(\Phi,\Phi)^{1/2}$ is a seminorm. For fixed $\xx$, we apply $B$ to the functions $\yy\mapsto K(\xx,\yy)$; for instance, $q^{{\rm st},L}_{3,3}(\xx,\xx)=B\left(K^{{\rm st},L}_{3,3}(\xx,\cdot),K^{{\rm st},L}_{3,3}(\xx,\cdot)\right)$ and by
\eqref{K-dec}, $q^{{\rm st},L}_{3,3}(\xx,\xx)=\left\|K^{{\rm hs},L}(\xx,\cdot)-H^L(\xx,\cdot)\right\|_B^2$. Hence
\begin{equation}\label{eq:comparison}
q^{{\rm st},L}_{3,3}(\xx,\xx)\le 2\,q^{{\rm hs},L}_{3,3}(\xx,\xx)+2\,\|H^L(\xx,\cdot)\|_B^2 .
\end{equation}
In order to estimate the corrector, and showing that $\|H^L(\xx,\cdot)\|_B^2$ will produce upper bounds that are not larger than the ones produced by $q^{{\rm hs},L}_{3,3}(\xx,\xx)$, we will employ an approach similar to the one in Otto\&Seis \cite[c.f. Proof of Proposition 2]{otto2011rayleigh}. We will first derive deterministic estimates for the corrector applied to an arbitrary function in $L^2(\D^{{\rm st}, L})$, then we will establish the desired bound at the level of the covariance $q_S^{{\rm st}, L}$.

\paragraph{The deterministic problem.} This part is inspired by  Otto\&Seis \cite[Proof of Proposition 2]{otto2011rayleigh}.
For $f\in L^2(\D^{{\rm st}, L})$, extended by zero to $\D^{{\rm hs}, L}$, let $w_f$ and $\widetilde{w}_f$ be the solutions to 
\begin{align*}
    &\Delta^2 w_f=\Ra \Delta' f, \qquad w_f=\partial_3 w_f=0 \mbox{ on } \{x_3=0\}\cup \{x_3=1\}\\
     &\Delta^2 \widetilde{w}_f=\Ra \Delta' f, \qquad \widetilde w_f=\partial_3 \widetilde w_f=0 \mbox{ on } \{x_3=0\}\text{ and decay as }x_3\rightarrow \infty\,.
\end{align*}
The corrector field
\[\hat{w}_f:=\widetilde{w}_f-w_f\quad \mbox{ on } \D^{{\rm st}, L}\]
solves
\[\Delta^2 \hat{w}_f=0 \quad  \hat{w}_f=\partial_3 \hat{w}_f=0 \mbox{ on } \{x_3=0\}\,.\]
We will see that values at the upper plate, $\hat{w}_f=\widetilde w_f$ and $\partial_3\hat{w}_f=\partial_3\widetilde w_f$ on $\{x_3=1\}$,
will not be needed.

Notice that 
\begin{equation}\label{rep:Htilde}
    \hat{w}_f(\xx)=\Ra\int_{\D^{{\rm st}, L}}H^L(\xx,\yy)f(\yy)\, {\rm d}\yy\,,
\end{equation}
where $H^L$ is defined by \eqref{HL}.

First observe by classical energy estimates, we have\footnote{Notice that the energy estimate for $\widetilde{w}$ is performed on $\D^{{\rm hs}, L}$!}
\begin{align*}
&\int_{\D^{{\rm st}, L}}|\Delta w_f|^2 {\rm d}\xx\leq \Ra^2 \int_{\D^{{\rm st}, L}}|f|^2\,{\rm d}\xx\\
&\int_{\D^{{\rm hs}, L}}|\Delta \widetilde{w}_f|^2 {\rm d}\xx\leq \Ra^2 \int_{\D^{{\rm hs}, L}}|f|^2 {\rm d}\xx=\Ra^2 \int_{\D^{{\rm st}, L}}|f|^2\,{\rm d}\xx.
\end{align*}
Since
\begin{align*}
    \int_{\D^{{\rm st}, L}}|\Delta \hat{w}_f|^2 {\rm d}\xx &\leq 2\int_{\D^{{\rm st}, L}}|\Delta \widetilde{w}_f|^2 {\rm d}\xx +2\int_{\D^{{\rm st}, L}}|\Delta w_f|^2 {\rm d}\xx\\
    &\leq 2\int_{\D^{{\rm hs}, L}}|\Delta \widetilde{w}_f|^2 {\rm d}\xx +2\int_{\D^{{\rm st}, L}}|\Delta w_f|^2{\rm d}\xx ,
\end{align*}
we have
\begin{equation}\label{Energy-estimate-corr}
\int_0^1\left\langle|\Delta \hat{w}_f|^2\right\rangle'\, {\rm d}x_3\leq 4\Ra^2 \int_0^1\left\langle|f|^2\right\rangle'\, {\rm d}x_3\,.
\end{equation}
\begin{lemma}[Deterministic corrector estimate]
    There exists a constant $C$, independent of $L>0$, such that, for every $f\in L^2(\D^{{\rm st}, L})$ it holds
    \begin{equation}\label{deterministic-estimate}
        \int_{\TT_L^{2}}\hat{w}_f(\xx',x_3)^2\, {\rm d}\xx'\leq C\Ra^2 x_3^4\|f\|^2_{L^2(\D^{{\rm st}, L})}\,, \qquad 0<x_3\leq \frac 12\,,
    \end{equation}
    or, equivalently,
    \[\left\langle\hat{w}_f^2\right\rangle'(x_3)\leq C \Ra^2 x_3^4\int_0^1\left\langle |f|^2\right\rangle' {\rm d}x_3,\qquad 0<x_3\leq \frac 12\,.\]
\end{lemma}
\begin{proof}
 Fist we observe that, since $\hat w_f=\partial_3\hat w_f=0$ on $x_3=0$, then, for any $x_3\in (0,1]$\footnote{These estimates come from the observation that, using the boundary conditions, $\hat{w}_f(\xx',x_3)=\int_0^{x_3}\int_0^y\partial_3^2 \hat{w}_f(\xx',z)\, {\rm d}z\, {\rm d}y=\int_0^{x_3}\int_z^{x_3}\partial_3^2 \hat{w}_f(\xx',z)\, {\rm d}y\, {\rm d}z=\int_0^{x_3}(x_3-z)\partial_3^2\hat{w}_f(\xx',z)\, {\rm d}z$. In particular
\[\langle \hat{w}_f\rangle'(x_3)\leq \int_0^{x_3}(x_3-z)\langle(\partial_{x_3}^2\hat{w}_f)^2\rangle'(z)\, {\rm d}z\,.\]
}
\begin{align}
&\left\langle \hat{w}_f^2\right\rangle'(x_3)\leq \frac{x_3^4}{4}\sup_{0<z<x_3}\left\langle \left(\partial_3^2 \hat{w}_f\right)^2\right\rangle'(z)\label{est:C1}\\
&\left\langle \hat{w}_f^2\right\rangle'(x_3)\leq \frac{x_3^3}{3}\int_0^{x_3}\left\langle \left(\partial_3^2 \hat{w}_f\right)^2\right\rangle'(z)\,{\rm d}z\,.\label{est:C2}
\end{align}
Also, the following Caccioppoli-type estimate holds
\begin{equation}\label{Caccioppoli}
    \sup_{0<x_3<\frac 18}\left\langle|\Delta \hat{w}_f|^2\right\rangle'(x_3)\leq C\int_0^1\left\langle|\nabla \hat{w}_f|^2\right\rangle'(z)\, {\rm d}z\,.
\end{equation}
For the proof of \eqref{Caccioppoli} we refer to Otto\&Seis \cite[Lemma 5]{otto2011rayleigh}.

\textbf{Case $0<x_3<\frac 18$}.
Then 
\begin{align*}
    \left\langle\hat{w}_f^2\right\rangle'(x_3)&\overset{\eqref{est:C1}}{\leq} \frac{x_3^4}{4}\sup_{0<x_3<\frac 18}\left\langle |\Delta\hat{w}_f|^2\right\rangle'(x_3)\\
    &\overset{\eqref{Caccioppoli}}{\leq} Cx_3^4\int_0^1\left\langle |\nabla \hat{w}_f|^2\right\rangle'(z) {\rm d}z\\
    &\leq C x_3^4\int_0^1\left\langle |\Delta\hat{w}_f|^2\right\rangle' (z){\rm d}z\leq C\Ra^2 x_3^4\int_0^1\left\langle |f|^2\right\rangle'(z){\rm d}z\,,
\end{align*}
where in the second-to-last inequality, a Poincar\'e estimate was used on each component of $\nabla \hat{w}_f$.

\textbf{Case $\frac 18\leq x_3 \leq\frac 12$}.
Here $x_3\geq 8^{-1}$, hence we have 
\begin{align*}
    \left\langle\hat{w}_f^2\right\rangle'(x_3)&\overset{\eqref{est:C2}}{\leq} \frac{x_3^3}{3}\int_0^1\left\langle|\Delta\hat{w}_f|^2\right\rangle'(z){\rm d}z\\
    &\overset{\eqref{Energy-estimate-corr}}{\leq} C\Ra^2x_3^3\int_0^1\left\langle |f|^2\right\rangle'(z){\rm d}z\leq 8Cx_3^4\Ra^2\int_0^1\langle |f|^2\rangle'(z){\rm d}z\,.
\end{align*}
\end{proof}
\paragraph{From the deterministic estimate to the covariance.}
Recall that $q_S^{{\rm st}, L}$ is the kernel associated to the integral operator $\mathcal{Q}_S$ on $L^2(\D^{{\rm st}, L})$ defined by 
\[(\mathcal{Q}_Sf)(\yy)=\int_{\D^{{\rm st},L}} q_S^{{\rm st}, L}(\yy,\zz)f(\zz)\, {\rm d}\zz\,.\] 
\begin{lemma}[Stochastic corrector estimate]\label{lem:stoch-corr}
    Let $q_S^{{\rm st}, L}$ be symmetric, positive semidefinite, horizontally stationary and satisfy
    \begin{equation}\label{qS:reduced-form}
         \left|q_S^{{\rm st}, L}(\yy,\zz)\right|\leq b(y_3)b(z_3)\qquad \mbox{ for all } \yy,\zz\in \D^{{\rm st}, L}
    \end{equation}
    for some $b\in L^2(0,1)$, $b\geq 0$. Then 
    \begin{equation}\label{est:stoc-error}
        \left\|H^L(\xx,\cdot)\right\|_{B}^2\leq C\Ra^2 x_3^4 Q_0,\qquad  0<x_3\leq \frac 12\,,
    \end{equation}
    where $Q_0=\int_0^1 b^2\, ds$ and $C$ is independent of $L$.
\end{lemma}
\begin{proof}
Recall that \[\left\|H^L(\xx,\cdot)\right\|_{B}^2=\Ra^2\int_{\D^{{\rm st},L}}\int_{\D^{{\rm st},L}}H^L(\xx,\yy)q_S^{{\rm st}, L}(\yy,\zz)H^L(\xx,\zz)\, {\rm d}\yy\,{\rm d}\zz \,.\]

Fix $x_3\in(0,\frac12]$, and for $\xx'\in\TT^2_L$ set
\[
r_{\xx'}:=\Ra\,H^L\big((\xx',x_3),\cdot\big),\qquad
F(\xx'):=\left\|H^L((\xx',x_3),\cdot)\right\|_B^2 .
\]
By \eqref{rep:Htilde}, $r_{\xx'}$ is the kernel of the map $f\mapsto\hat w_f(\xx',x_3)$, that is,
$\hat w_f(\xx',x_3)=\langle r_{\xx'},f\rangle$. Moreover, $r_{\xx'}\in L^2(\D^{{\rm st},L})$ by Lemma~\ref{lem:representer}:
since the energy estimate implies that $f\mapsto\hat w_f(\xx',x_3)$ is a bounded functional on $L^2(\D^{{\rm st},L})$, its
kernel is square integrable. Hence, for all $f\in L^2(\D^{{\rm st},L})$,
\begin{equation}\label{eq:repr-summary}
\hat w_f(\xx',x_3)=\langle r_{\xx'},f\rangle,\qquad F(\xx')=\langle r_{\xx'},\mathcal{Q}_Sr_{\xx'}\rangle,
\end{equation}
where $\mathcal{Q}_S$ is the operator on $L^2(\D^{{\rm st},L})$ with kernel $q^{{\rm st},L}_S$ and $\langle\cdot,\cdot\rangle$ the inner product of $L^2(\D^{{\rm st},L})$.

\medskip

\emph{Step 1: Expansion in the eigenfunctions of $\mathcal{Q}_S$.} By \eqref{qS:reduced-form}, the kernel $q^{{\rm st},L}_S$ belongs to $L^2(\D^{{\rm st},L}\times \D^{{\rm st},L})$, hence $\mathcal{Q}_S$ is a compact operator\footnote{An integral operator on $L^2$ whose kernel lies in $L^2$ of the product space is a Hilbert–Schmidt operator, and Hilbert–Schmidt operators are compact.}; it is self-adjoint since $q^{{\rm st},L}_S$ is symmetric. Hence there is an orthonormal basis $(\varphi_j)_j$ of $L^2(\D^{{\rm st},L})$ with $\mathcal{Q}_S\varphi_j=\mu_j\varphi_j$, and $\mu_j=\langle\varphi_j,\mathcal{Q}_S\varphi_j\rangle\ge0$ since $q^{{\rm st},L}_S$ is positive semidefinite. Expanding $r_{\xx'}$ in this basis and using \eqref{eq:repr-summary} with
$f=\varphi_j$, we obtain
\begin{align}\label{eq:F-expansion}
F(\xx')&=\langle r_{\xx'}, \mathcal{Q}_S r_{\xx'}\rangle \notag\\
&=\langle r_{\xx'}, \sum_{j}\langle r_{\xx'},\varphi_j\rangle \mathcal{Q}_S\varphi_j\rangle= \langle r_{\xx'}, \sum_{j}\langle r_{\xx'},\varphi_j\rangle \mu_j\varphi_j\rangle\notag\\
&=\sum_{j}\langle r_{\xx'}, \varphi_j\rangle \mu_j\langle r_{\xx'},\varphi_j\rangle\notag\\
&=\sum_{j} \mu_j\langle r_{\xx'},\varphi_j\rangle^2.
\end{align}

\emph{Step 2: The averaged bound.} The deterministic estimate \eqref{deterministic-estimate}, multiplied by
$L^2$, states that $\int_{\TT^2_L}\hat w_f(\xx',x_3)^2\,{\rm d}\xx'\le C\,\Ra^2x_3^4\,\|f\|^2_{L^2(\D^{{\rm st},L})}$ for
every $f\in L^2(\D^{{\rm st},L})$. Applied to $f=\varphi_j$, with $\|\varphi_j\|_{L^2}=1$, and summed according to
\eqref{eq:F-expansion}, where all terms are nonnegative, it gives
\[
\int_{\TT^2_L}F(\xx')\,{\rm d}\xx'=\sum_j\lambda_j\int_{\TT^2_L}\hat w_{\varphi_j}(\xx',x_3)^2\,{\rm d}\xx'
\le C\,\Ra^2x_3^4\sum_j\lambda_j\le C\,\Ra^2x_3^4\,L^2Q_0,
\]
where in the last step we used Lemma~\ref{lem:eigen-sum}.

\emph{Step 3: From the average to every point.} Since $q^{{\rm st},L}_S$ is horizontally stationary and
$H^L$ is invariant under simultaneous horizontal translations of both variables, $F$ does not
depend on $\xx'$. Hence $F(\xx')=L^{-2}\int_{\TT^2_L}F$ for every $\xx'$, and Step 2 gives
\eqref{est:stoc-error}.
\end{proof}

Finally we are in the position to prove Proposition \ref{Pr1}:
\begin{proof}[Proof of Proposition \ref{Pr1}]
    By Lemma \ref{lem:stoch-corr} and \eqref{eq:comparison} we have 
    \begin{equation*}
        q_{3,3}^{{\rm st}, L}(\xx,\xx)\leq 2  q_{3,3}^{{\rm hs}, L}(\xx,\xx)+2 CQ_0 \Ra^2 x_3^4\,.
    \end{equation*}
    Notice that the assumptions in Case 1 and Case 2 in Proposition \ref{Pr1} can be reduced to the form \eqref{qS:reduced-form}. In fact, $\sup_{\rr'}\sum_{\kk'}(1+|\rr'+L\kk'|/\ell)^{-N}\leq C_N$ and all the scales $\ell$ appearing in $(1, \sqrt{y_3z_3}, \sqrt{(1-y_3)(1-z_3)})$ are at most of the size $1$\footnote{Indeed, choose $\kk'_0\in\mathbb{Z}^2$ with $|\rr'+L\kk'_0|\le L/\sqrt2$. For $\kk'\neq\kk_0$, the
reverse triangle inequality gives $|\rr'+L\kk'|\ge L|\kk'-\kk'_0|-L/\sqrt2\ge cL|\kk'-\kk'_0|$ with
$c:=1-\frac1{\sqrt2}$, hence
\[
\sum_{\kk'}\left(1+\frac{|\rr'+L\kk'|}{\ell}\right)^{-N}\leq 1+\left(\frac{\ell}{cL}\right)^N\sum_{\kk'\neq\kk'_0}|\kk'-\kk'_0|^{-N}
\overset{\ell\leq 1\leq  L}{\leq}1+c^{-N}\sum_{\kk'\in\mathbb{Z}^2\setminus\{0\}}|\kk'|^{-N}=:C_N,
\]
which is finite since $N>2$.}. Hence we see that the corrector terms do not alter the estimates obtained in the periodic half-space:
\begin{itemize}
    \item In Case 1, $b(s)=C(s(1-s))^{\beta}$, hence $Q_0\lesssim \int_0^1 (s(1-s))^{2\beta}\, {\rm d}s\lesssim 1$. Combining the estimates from Lemma \ref{lemma:bulkCov}, Lemma \ref{lemma:BL-bulk}, Lemma \ref{extension-lemma} and \eqref{eq:comparison}, the final estimate is
\begin{equation*}
    \left|q_{3,3}^{{\rm st},L}(\xx,\xx)\right| \lesssim \Ra^2 x_3^4.
\end{equation*}
\item In Case 2, the three profiles can again be reformulated into \eqref{qS:reduced-form}, with 
$$b^{\rm bot}(s)=C s^{\beta} \chi_{s\leq \eta}, \quad b^{\rm top}(s)=C (1-s)^{\beta} \chi_{s\geq  1-\eta} \quad b^{\rm bulk}(s)=C (s(1-s))^{1+\sigma} \chi_{\eta<s< 1-\eta}\,,$$
such that 
$$Q_0\lesssim \eta^{2\beta+1}+\eta^{2\beta}+\eta^{2\beta+1}\lesssim \eta^{2\beta}\,,$$
since $\eta\ll 1$ by assumption.
Again combining the estimates from Lemma \ref{lemma:bulkCov}, Lemma \ref{lemma:BL-bulk}, Lemma \ref{extension-lemma} and \eqref{eq:comparison}, the final estimates is
\begin{equation*}
    \left|q_{3,3}^{{\rm str},L}(\xx,\xx)\right| \lesssim \Ra^2\left(x_3^{2\beta + 4} + \eta^{2\beta} x_3^4\right).
\end{equation*}
\end{itemize}
\end{proof}

\subsection{Proof of Theorem \ref{Main Theorem} }\label{sec:proofMT}
\begin{proof}[Proof of Theorem \ref{Main Theorem}]
   Recall the upper bound \eqref{Nu-apriori}, namely
    \begin{equation}\label{Nu_UB}	
		\Nu \lesssim  \frac{1}{\delta} + \frac{1}{\delta^2}\int_{0}^{\delta}  |\widetilde{q}_{3,3}(x_3)|^{2} \, {\rm d}x_3\,.
	\end{equation}
    Recalling that $\Ra\gg 1$, we will now apply Proposition \ref{Pr1} with $\beta=1$.
\begin{description}
    \item[Case 1:] 
    Inserting \eqref{q33-bulk} into \eqref{Nu_UB} we obtain
\begin{equation*}
	\Nu \lesssim \frac{1}{\delta} + \frac{1}{\delta^2}\int_{0}^{\delta} \Ra^{4} x_3^8\, {\rm d}x_3.
\end{equation*}
The optimal $\delta$ is 
$$\delta :=\Ra^{-\frac 12}$$
yielding
$$\Nu \lesssim \Ra^{\frac{1}{2}}\,.$$

\item[Case 2: ]
As in Case 1, we now insert \eqref{q33-BL} into \eqref{Nu_UB}. Choosing $\beta=1$, for every $0< \delta\leq \frac 12$, we have
\begin{equation}\label{opt}
	\Nu \lesssim\frac{1}{\delta} + \frac{1}{\delta^2}\int_{0}^{\delta} \Ra^{4} \left(x_3^{6}+\eta^{2}x_3^4 \right)^{2}\, {\rm d}x_3\\
    \leq \frac{1}{\delta} +\Ra^4\delta^7 \left(\delta^{4}+\eta^{4}\right)=\frac{1}{\delta}+\Ra^{4}\delta^{11}+\Ra^4 \eta^4\delta^7.
\end{equation}
Optimizing in $\delta$ we choose 
\[
\delta:=\min\left\{\Ra^{-1/3},\,(\Ra\,\eta)^{-1/2}\right\}\leq \frac 12, \qquad \mbox{ for } \Ra\geq 8
\]
Inserting in \eqref{opt} we obtain 
\[\Nu\lesssim \max\left\{\Ra^{\frac 13}, (\eta\Ra)^{\frac 12}\right\}\,,\]
implying 
\[\Nu \lesssim \Ra^{\frac 13}\,,\]
if 
$\eta\leq \Ra^{-\frac 13}$.
\end{description}
\end{proof} 

\appendix
\section{Appendix}

\subsection{Nusselt number estimate}\label{app:nusselt}

\begin{lemma}
	Let $\delta \in (0,1]$ be arbitrary, then the Nusselt number \eqref{eq:stochasticNu} satisfies
	\begin{equation*}
		\Nu \lesssim \frac{1}{\delta} + \frac{1}{\delta^2}\sum_{j=1}^3 \int_0^{\delta} \left\langle |\widetilde{q}_{3,3}(x_3)|^2\right\rangle'{\rm d}x_3.
	\end{equation*}
\end{lemma}
\begin{proof}
For $T>0$ we denote by
\begin{equation*}
  \langle f\rangle_T(x_3) := \frac{1}{T}\int_0^T \frac{1}{L^2}\int_{[0,L)^2} f(s,\boldsymbol{x}',x_3)\,\mathrm{d}\boldsymbol{x}'\,\mathrm{d}s 
\end{equation*}
the finite-time horizontal average.
Applied to \eqref{eq:EvolExp}, we obtain
	\begin{equation*}
		 \partial_{3} \left\langle \left[\left({\rm Id}+\frac 12 \widetilde{Q}_{\ub}\right)\nabla E\right]^{(3)}\right\rangle_T=\frac{1}{TL^2}\int_{[0,L)^2}E(t,\xx',x_3)-E(0,\xx',x_3)\, {\rm d}\xx'\,
	\end{equation*}    
    with $[\cdot]^{(3)}$ denoting the third component of the vector. Denote
    \[R_T(x_3):=\frac{1}{TL^2}\int_{[0,L)^2}E(t,\xx',x_3)-E(0,\xx',x_3)\, {\rm d}\xx',\]
    then by the maximum principle we observe that
    \begin{equation}\label{eq:est-RT}
        |R_T|\leq \frac{1}{T}.
    \end{equation}
	Integrating the last identity in $x_3$ we have 
    	\begin{equation*}
		 \left\langle \left[\left({\rm Id}+\frac 12 \widetilde{Q}_{\ub}\right)\nabla E\right]^{(3)}\right\rangle_T (x_3) = \left.\left\langle \left[\left({\rm Id}+\frac 12 \widetilde{Q}_{\ub}\right)\nabla E\right]^{(3)}\right\rangle\right|_{x_3=0}+R_T(x_3)\,.
	\end{equation*}  
    Due to the boundary conditions on $\widetilde{Q}_\ub$ we have
	\[\left.\left\langle \left[\left({\rm Id}+\frac 12 \widetilde{Q}_{\ub}\right)\nabla E\right]^{(3)}\right\rangle_T\right|_{x_3=0}=\left\langle\partial_3 E\right\rangle_T|_{x_3=0}=: - \Nu_T\,\]
    and, by \eqref{eq:stochasticNu}, it holds  
    \[\Nu=\lim_{T\rightarrow \infty}\Nu_T\,.\]
	Therefore for any $\delta \in(0,1]$ it holds
	\begin{equation*}
		- \Nu_T = \frac{1}{\delta} \int_0^{\delta} \left\langle \left[\left({\rm Id}+\frac 12 \widetilde{Q}_{\ub}\right)\nabla E\right]^{(3)}\right\rangle_T(x_3) - R_T(x_3) \,{\rm d}x_3\,.
	\end{equation*}
	Next using the horizontal stationarity of the coefficients $\widetilde{q}_{i,j}$, which depend on $x_3$ only, we obtain
    \begin{align*} 
    &\left\langle \left[\left({\rm Id}+\frac 12 \widetilde{Q}_{\ub}\right)\nabla E\right]^{(3)}\right\rangle_T(x_3)\\
    &=\frac 12\widetilde{q}_{3,1}(x_3)  \langle  \partial_1 E\rangle_T(x_3)+\frac 12\widetilde{q}_{3,2}(x_3)\langle   \partial_2 E\rangle_T(x_3)+\left(1+\frac 12 \widetilde{q}_{3,3}(x_3)\right)\langle  \partial_3 E\rangle_T(x_3)\\
    &=\left(1+\frac 12 \widetilde{q}_{3,3}(x_3)\right)\langle  \partial_3 E\rangle_T(x_3)\,,
    \end{align*}
    where we used that $\langle \partial_i E\rangle_T=0$ for $i=1,2$ due to periodicity.
    Hence, the localized Nusselt number identity reduces to 
    	\begin{equation*}
		- \Nu_T = \frac{1}{\delta} \int_0^{\delta} \left(1+\frac 12 \widetilde{q}_{3,3}(x_3)\right)\langle  \partial_3 E\rangle_T(x_3) - R_T(x_3) \,{\rm d}x_3.
	\end{equation*}
    Now observe that, on the one hand
    \[\int_0^{\delta}\langle \partial_3 E\rangle_T(x_3)\, {\rm d}x_3=\langle E\rangle_T(\delta)-\langle E\rangle_T(0)\geq -1\,,\]
    because of the boundary conditions and the fact that $E\geq 0$ by the maximum principle.
     On the other hand, by the Cauchy--Schwarz inequality in $x_3$ and Jensen's inequality for the normalized space-time average, $|\langle\partial_3E\rangle_T|^2\le \langle|\partial_3E|^2\rangle_T$, 
    \[\int_0^{\delta}\widetilde{q}_{3,3}(x_3)\langle\partial_3E\rangle_T (x_3){\rm d}x_3\leq \left(\int_0^{\delta}\left|\frac 12\widetilde{q}_{3,3}(x_3)\right|^2\, {\rm d}x_3\right)^{\frac 12}\left(\int_0^{\delta}\langle|\partial_3E|^2\rangle_T(x_3) \,{\rm d}x_3\right)^{\frac 12}\,.\]
    Since, by the energy identity at fixed $T$ (c.f. proof of \cite[Lemma 3.2]{FLN2026a})
    \[\int_0^{\delta}\langle | \partial_3 E|^2 \rangle_T(x_3) {\rm d}x_3\leq \int_0^1 \left \langle \left| \left( {\rm Id} + \frac{1}{2}\widetilde{Q}_\ub\right)^{\frac{1}{2}}\nabla E\right|^2\right \rangle_T(x_3){\rm d}x_3\leq \Nu_T+\frac{1}{T},\]
	we conclude with \eqref{eq:est-RT} that
	\[-\Nu_T\geq -\frac{1}{\delta}-\frac 12\Nu_T-\frac 18 \frac{1}{\delta^2}\int_0^{\delta}|\widetilde{q}_{3,3}(x_3)|^2\, {\rm d}x_3-\frac{3}{2T}\,,\]
    implying 
    \[\Nu_T\lesssim\frac{1}{\delta}+\frac{1}{\delta^2}\int_0^{\delta}|\widetilde{q}_{3,3}(x_3)|^2\, {\rm d}x_3+\frac{3}{T} \]
    and, taking the $\limsup$ as $T\rightarrow \infty$ we obtain \eqref{Nu-apriori}.
\end{proof}

\subsection{Green functions estimates}

\begin{lemma}[Far-field bounds for Green's function]\label{appLemma:GreenEst}
	The Green kernel $G$ defined in \eqref{green-kernel} satisfies the horizontal bounds
	$$|G(\xib,\etab)| \lesssim \frac{\xi_3^2 \eta_3^2}{|\xib'-\etab'|^3},
	\quad
	|\nabla_{\etab'}G(\xib,\etab)| \lesssim \frac{\xi_3^2 \eta_3^2}{|\xib'-\etab'|^4}, \quad \left|\partial_{\eta_3}G(\xib,\etab)\right| \lesssim \frac{\xi_3^2 \eta_3}{|\xib'-\etab'|^3}$$
	in the far-field regime $|\xib-\etab|>8\min\{\xi_3,\eta_3\}$.
\end{lemma}
\begin{proof}
	Set
	\begin{align*}
	    \rho &:=|\xib'-\etab'|,\\
        A &:=|\xib^*-\etab|=\sqrt{\rho^2+(\xi_3+\eta_3)^2},\\
        B &:=|\xib-\etab|=\sqrt{\rho^2+(\xi_3-\eta_3)^2}.
	\end{align*}
	Since
	$$A^2-B^2=(\xi_3+\eta_3)^2-(\xi_3-\eta_3)^2=4\xi_3\eta_3,$$
	it holds
	$$A-B= \frac{A^2 - B^2}{A+B} = \frac{4\xi_3\eta_3}{A+B},$$
	and therefore
	\begin{equation}\label{rep-G}
	G(\xib,\etab)=\frac{1}{16\pi}\frac{(A-B)^2}{A}=\frac{\xi_3^2\eta_3^2}{\pi\,A\,(A+B)^2}.
	\end{equation}
	In the far-field we use the rough lower bounds
	$$A\ge \rho,\qquad A+B\ge 2\rho.$$
	Then from the representation above,
	$$|G(\xib,\etab)|\lesssim\frac{\xi_3^2\eta_3^2}{\rho^3}.$$
	Similarly, differentiating
	$$G(\xib,\etab)=\frac{\xi_3^2\eta_3^2}{\pi\,A\,(A+B)^2},$$
	we obtain
	$$|\nabla_{\etab}'G|
	\lesssim\xi_3^2\eta_3^2\left(\frac{|\nabla_{\etab}'A|}{A^2(A+B)^2}+
	\frac{|\nabla_{\etab}'A|+|\nabla_{\etab}'B|}{A(A+B)^3}\right).$$
    Further note that
    \[|\nabla_{\etab}'A|\leq \frac{|\xib' - \etab'|}{A} \leq 1\]
    and similarly for $B$ it holds $|\nabla_{\etab}'B|\leq 1$. Together with the crude lower bounds above we obtain
	$$|\nabla_{\etab}'G|\lesssim\frac{\xi_3^2\eta_3^2}{\rho^4}.$$
	Finally, we obtain
    \[\partial_{\eta_3}G(\xib, \etab) = \frac{2\xi_3^2 \eta_3}{\pi A(A+B)^2} + \frac{\xi_3^2\eta_3^2}{\pi}\partial_{\eta_3}\frac{1}{A(A+B)^2}\]
    where one can easily check that the first summand dominates. Thus we conclude as before
    \[\left|\partial_{\eta_3}G(\xib, \etab)\right|\leq \frac{2\xi_3^2 \eta_3}{\pi \rho^3}.\]
\end{proof}
The following result, already stated in Section \ref{preliminaries} contains near- and far-field upper bounds for $K_{3,3}$ that will be crucially used in Lemma \ref{lemma:BL-bulk}, and Lemma \ref{lemma:bulkCov}.
\begin{lemma}[Near-/far-field kernel bounds]\label{appendix:KEstimates}
	For $\xi_3, \eta_3>0$, $\xib\neq \etab$ the kernel $K_{3,3}$ satisfies:
    \begin{itemize}
        \item near-field $|\xib-\etab|\leq 8\min\{\xi_3,\eta_3\}$:
            \begin{equation*}
                |K_{3,3}(\xib,\etab)|\leq \frac{C}{|\xib-\etab|};
            \end{equation*}
        \item far-field $|\xib-\etab|>8\min\{\xi_3,\eta_3\}$:
            \begin{equation*}
                |K_{3,3}(\xib,\etab)|\leq \frac{C\min\{\xi_3,\eta_3\}^2}{|\xib-\etab|^3}.
            \end{equation*}
    \end{itemize}
\end{lemma}

\begin{proof}
	The proof of the near-field estimate can be found in \cite[Prop. 2.6]{ChanMal2025}. For the bound in the far-field $|\xib-\etab|\geq 8\min\{\xi_3,\eta_3\}$ we start again from the representation \eqref{rep-G}, i.e.
    $$G(\xib,\etab)=\frac{\xi_3^2\eta_3^2}{\pi A(A+B)^2}\qquad \mbox{ with } A=|\xib^{\star}-\etab|, \quad B=|\xib-\etab|\,.$$
It is clear that $B\leq A$. W.l.o.g. we can assume $\xi_3\leq \eta_3$, then
$$\xi_3+\eta_3=(\eta_3-\xi_3)+2\xi_3=|\xi_3-\eta_3|+2\min\{\xi_3,\eta_3\}\,.$$
By the triangle inequality we have 
$$A \leq \sqrt{\rho^2 + (\xi_3 - \eta_3)^2} + 2|\min\{\xi_3,\eta_3\}| = B + 2\min\{\xi_3,\eta_3\} \leq \frac{5}{4}B\,.$$
Hence, we deduce that, in the far field $B\sim A$.

Now we focus on the horizontal part of $G$. Consider
$$F(\rho)=\frac{1}{A(A+B)^2}$$
with $\rho=|\xib'-\etab'|$, then with our previous discussion observe that it holds
\[\Delta_{\etab}'F = F''(\rho) + \frac{F'(\rho)}{\rho}\]
and
\[|F(\rho)| \sim B^{-3},\quad \frac{|F'(\rho)|}{\rho} \lesssim B^{-5},\quad|F''(\rho)| \lesssim B^{-5}.\]
With this we deduce that $K_{3,3} = \Delta_{\etab'}G$ satisfies the upper bound
    \[|K_{3,3}(\xib,\etab)| \lesssim \frac{\xi_3^2\eta_3^2}{|\xib-\etab|^5}.\]
    Now we observe that
    \[\xi_3^2 \eta_3^2 = \min(\xi_3,\eta_3)^2 \max(\xi_3,\eta_3)^2.\]
    First of all, assume $\xi_3 \leq \eta_3$. If $\eta_3 - \xi_3 \geq \frac{1}{2}\eta_3$, then
    \[|\xib - \etab|^2 = |\xib' - \etab'|^2 + (\xi_3 - \eta_3)^2 \gtrsim \eta_3^2.\]
    If $\eta_3 - \xi_3 \leq \frac{1}{2}\eta_3$, then we use the far-field estimate
    \[|\xib - \etab|^2 \geq (8 \min(\xi_3,\eta_3))^2 \gtrsim \eta_3^2.\]
    Analogously, in the case $\eta_3 \leq \xi_3$, we obtain $\xi_3^2 \lesssim |\xib - \etab|^2$, hence we have shown
    \[\max(\xi_3,\eta_3)^2 \lesssim |\xib-\etab|^2\]
    and hence
    \[\xi_3^2\eta_3^2 \lesssim \min(\xi_3,\eta_3)^2|\xib - \etab|^2.\]
    In conclusion, it holds
    \[|K_{3,3}(\xib,\etab)| \lesssim \frac{\min(\xi_3,\eta_3)^2}{|\xib-\etab|^3}.\]
\end{proof}

\section{Two auxiliary lemmas for the corrector estimate}

\begin{lemma}[The corrector kernel as a representer]\label{lem:representer}
Let $q^{{\rm st},L}_S$ be symmetric, positive semidefinite, and satisfy \eqref{qS:reduced-form}, and let $\mathcal{Q}_S$ be the
operator on $L^2(\D^{{\rm st},L})$ with kernel $q^{{\rm st},L}_S$. Then, for every $\xx\in \D^{{\rm st},L}$, the
function $r_{\xx}:=\Ra\, H^L(\xx,\cdot)$ belongs to $L^2(\D^{{\rm st},L})$, and
\[
\hat w_f(\xx)=\langle r_{\xx},f\rangle\quad\text{for all } f\in L^2(\D^{{\rm st},L}),\qquad
\| H^L(\xx,\cdot)\|_B^2=\langle r_{\xx},\mathcal{Q}_S r_{\xx}\rangle .
\]
\end{lemma}

\begin{proof}
Fix $\xx\in \D^{{\rm st},L}$ and consider the linear map $\ell_{\xx}(f):=\hat w_f(\xx)$ on $L^2(\D^{{\rm st},L})$. Then first of all, observe that it holds $\hat w_f\in H^2(\D^{{\rm st},L})$, with norm controlled by $f$. Indeed: by \eqref{Energy-estimate-corr}
\[\|\Delta\hat w_f\|_{L^2}\le2\,\Ra\,\|f\|_{L^2}.\]
Since $\hat w_f=\partial_3\hat w_f=0$ on $\{x_3=0\}$, also the
horizontal derivatives of $\hat w_f$ vanish there, hence
\[\nabla\hat w_f=0 \text{ on }\{x_3=0\}.\]
Applying the Poincar\'e inequality  first to each component of $\nabla\hat w_f$ and then to $\hat w_f$ gives
\[\|\nabla\hat w_f\|^2_{L^2}\lesssim\|\nabla^2\hat w_f\|^2_{L^2}\]
and
\[\|\hat w_f\|^2_{L^2}\lesssim \|\nabla^2\hat w_f\|^2_{L^2}.\]
Hence
\[\|\hat w_f\|_{H^2}\le C\,\Ra\,\|f\|_{L^2}.\]
By the Sobolev embedding $H^2(\D^{{\rm st},L})\hookrightarrow C^0$ in three
dimensions,
\[|\ell_{\xx}(f)|=|\hat w_f(\xx)|\leq \sup_{\D^{{\rm st}, L}}|\hat{w}|\leq C_{\rm Sob} \|\hat{w}\|_{H^2(\D^{{\rm st}, L})}\leq  C_{\rm Sob}\,C\,\Ra\,\|f\|_{L^2}<\infty,\]
hence $\ell_{\xx}$ is a bounded linear functional on $L^2(\D^{{\rm st},L})$. By the Riesz representation theorem, there exists a unique $r\in L^2(\D^{{\rm st},L})$ with 
\[\ell_{\xx}(f)=\langle r,f\rangle \text{ for all }f \in L^2(\D^{{\rm st},L}).\]
By \eqref{rep:Htilde}, it furthermore holds
\[\ell_{\xx}(f)=\int_{\D^{{\rm st},L}} \Ra\, H^L(\xx,\yy)f(\yy)\,{\rm d}\yy,\]
hence
\[\int_{\D^{{\rm st},L}} (r-\Ra\,H^L(\xx,\cdot))f=0\]
for all bounded $f$, which implies $r= r_{\xx}$ a.e. In particular, $ H^L(\xx,\cdot)\in L^2(\D^{{\rm st},L})$. By \eqref{qS:reduced-form} and $r_{\xx}\in L^2(\D^{{\rm st},L})$, the double integral defining $\| H^L(\xx,\cdot)\|_B^2$ converges absolutely, and 
\[\| H^L(\xx,\cdot)\|_B^2=\int_{\D^{{\rm st},L}} \int_{\D^{{\rm st},L}} r_{\xx}(\yy)\,q^{{\rm st},L}_S(\yy,\zz)\,r_{\xx}(\zz)\,{\rm d}\yy\,{\rm d}\zz=\langle 
r_{\xx},\mathcal{Q}_S r_{\xx}\rangle.\]
\end{proof}

\begin{lemma}[Sum of the eigenvalues]\label{lem:eigen-sum}
Under the assumptions of Lemma~\ref{lem:representer}, let $(\mu_j)_j$ be the eigenvalues of $\mathcal{Q}_S$, counted with multiplicity. Then $\sum_j\mu_j\le L^2Q_0$, with $Q_0=\int_0^1b^2\,{\rm d}s$.
\end{lemma}

\begin{proof}
Let $(\varphi_j)_j$ be an orthonormal eigenbasis of $\mathcal{Q}_S$ with eigenvalues $(\mu_j)_j$.
We partition $\D^{{\rm st},L}$ into $N_n$ boxes $c_n$, of size $L2^{-n}\times L2^{-n}\times2^{-n}$. Next set $g_n:=|c_n|^{-1/2}\chi_{c_n}$ such that
\[\langle g_n, \mathcal{Q}_S g_n\rangle=\frac{1}{|c_n|}\int_{c_n}\int_{c_n} q^{{\rm st}, L}_S(\yy,\zz){\rm d}\yy {\rm d}\zz\]
and let $P_n$ be the orthogonal projection onto the span of the $g_n$, i.e. $\|P_n \phi\|=\sum_{i=1}^{N_n} \langle g_i,\phi\rangle^2$. Expanding $g_n$ in the eigenbasis yields
\begin{align*}
    \sum_{i=1}^{N_n}\langle g_i,\mathcal{Q}_S g_i\rangle&=
\sum_{i=1}^{N_n}\langle g_i, \sum_j\langle g_i,\varphi_j\rangle \mathcal{Q}_S\varphi_j\rangle\\
&=
\sum_{i=1}^{N_n}\langle g_i, \sum_j\langle g_i,\varphi_j\rangle \mu_j\varphi_j\rangle\\
&=
\sum_{i=1}^{N_n}\sum_j\mu_j\langle g_i,\varphi_j\rangle^2=\sum_j\mu_j\|P_n\varphi_j\|^2.
\end{align*}
On the other hand, by \eqref{qS:reduced-form} and the Cauchy-Schwarz inequality on each box we obtain
\[
\langle g_n,\mathcal{Q}_S g_n\rangle=\frac1{|c_n|}\int_{c_n} \int_{c_n} q^{{\rm st},L}_S(\yy,\zz)\,{\rm d}\yy\,{\rm d}\zz\le\frac1{|c_n|}\left(\int_{c_n}b(y_3)\,{\rm d}\yy\right)^2
\le\int_{c_n} b(y_3)^2\,{\rm d}\yy ,
\]
thus, combining all together, we have  \[\sum_j\mu_j\|P_n\varphi_j\|^2\le L^2Q_0 \quad \forall \, n\,.\]
Since $P_n\phi\rightarrow\phi$ in $L^2(\D^{{\rm st},L})$ for every
$\phi$ (since box averages of a continuous function converge uniformly, and continuous functions are dense in $L^2(\D^{{\rm st}, L})$), it holds $\|P_n\varphi_j\|_{L^2}\rightarrow 1$, by and Fatou's lemma we deduce that
\[\sum_j\mu_j\le\liminf_n\sum_j\mu_j\|P_n\varphi_j\|^2\le L^2Q_0.\]
\end{proof}

\section*{Acknowledgements}
During this project, FF and TL have received funding from the European Research Council (ERC) under the EU-Horizon Europe ERC-2021-ADG research and innovation programme (project ``Noise in Fluids", grant agreement no. 101053472), while TL has received further funding from the Deutsche Forschungsgemeinschaft (DFG, German Research Foundation) - project number 569224218. CN gratefully acknowledges the support of TRR181 and the Scuola Normale Superiore (Pisa), which enabled mutual visits in 2021–2022, during which the project was initiated in collaboration with FF.

\section*{AI Usage Disclosure}
The model, the main ideas and the strategy of the proofs were developed by the authors, in a project started in 2021. During the preparation of the manuscript, the authors used Claude (Anthropic) and DeepSeek as an assistant, to check proofs and computations and to help work out some technical steps. In particular, following the authors' idea of adapting the approach of \cite{otto2011rayleigh} to the corrector estimate, Claude helped to develop the details of the energy-based proof in Section \ref{sec:corrector}, in particular Lemmas \ref{lem:representer} and \ref{lem:eigen-sum} in the appendix and it helped to adapt the half-space estimates, which the authors had derived for non-periodic covariances, to horizontally periodic ones. Claude was also used to improve the language and presentation of the introduction and the abstract. All suggestions were critically assessed and verified in full by the authors, who take full responsibility for the content of the publication.

\bibliographystyle{plain}
\bibliography{bibliography}

@article{K1962,
  title={Turbulent thermal convection at arbitrary {P}randtl number},
  author={Kraichnan, Robert H.},
  journal={The Physics of Fluids},
  volume={5},
  number={11},
  pages={1374--1389},
  year={1962},
  publisher={AIP Publishing}
}

@article{Sp1963,
  title={A Generalization of the Mixing-Length Theory of Turbulent Convection.},
  author={Spiegel, Edward A.},
  journal={Astrophysical Journal, vol. 138, p. 216},
  volume={138},
  pages={216},
  year={1963}
}

@article{Malkus1954,
  title={The heat transport and spectrum of thermal turbulence},
  author={Malkus, Willem V. R.},
  journal={Proceedings of the Royal Society of London. Series A. Mathematical and Physical Sciences},
  volume={225},
  number={1161},
  pages={196--212},
  year={1954},
  publisher={The Royal Society London}
}

@article{choffrutNobiliOttoUpperBoundsOnNusseltNumberAtFinitePrandtlNumber,
    title = {Upper bounds on {N}usselt number at finite {P}randtl number},
    journal = {Journal of Differential Equations},
    volume = {260},
    number = {4},
    pages = {3860-3880},
    year = {2016},
    IGNOREissn = {0022-0396},
    doi = {https://doi.org/10.1016/j.jde.2015.10.051},
    url = {https://www.sciencedirect.com/science/article/pii/S0022039615006051},
    author = {Antoine Choffrut and Camilla Nobili and Felix Otto}
}

@article{DC96,
  title={Variational bounds on energy dissipation in incompressible flows. {III}. {C}onvection},
  author={Doering, Charles R. and Constantin, Peter},
  journal={Physical Review E},
  volume={53},
  number={6},
  pages={5957},
  year={1996},
  publisher={APS}
}

@article{DOR06,
  title={Bounds on vertical heat transport for infinite-{P}randtl-number {R}ayleigh-{B}{\'e}nard convection},
  author={Doering, Charles R. and Otto, Felix and Reznikoff, Maria G.},
  journal={Journal of Fluid Mechanics},
  volume={560},
  pages={229--241},
  year={2006},
  publisher={Cambridge University Press}
}

@article{NO17,
  title={Limitations of the background field method applied to {R}ayleigh-{B}{\'e}nard convection},
  author={Nobili, Camilla and Otto, Felix},
  journal={Journal of Mathematical Physics},
  volume={58},
  number={9},
  pages={093102},
  year={2017},
  publisher={AIP Publishing LLC}
}

@article{CD99,
  title={Infinite {P}randtl number convection},
  author={Constantin, Peter and Doering, Charles R.},
  journal={Journal of Statistical Physics},
  volume={94},
  number={1},
  pages={159--172},
  year={1999},
  publisher={Springer}
}

@article{AGL09,
  title={Heat transfer and large scale dynamics in turbulent {R}ayleigh-{B}{\'e}nard convection},
  author={Ahlers, Guenter and Grossmann, Siegfried and Lohse, Detlef},
  journal={Reviews of Modern Physics},
  volume={81},
  number={2},
  pages={503},
  year={2009},
  publisher={APS}
}

@article{doering2019,
  title={Absence of evidence for the ultimate regime in two-dimensional {R}ayleigh-{B}{\'e}nard convection},
  author={Doering, Charles R. and Toppaladoddi, Srikanth and Wettlaufer, John S.},
  journal={Physical Review Letters},
  volume={123},
  number={25},
  pages={259401},
  year={2019},
  publisher={APS}
}

@article{nobili2023role,
  title={The role of boundary conditions in scaling laws for turbulent heat transport},
  author={Nobili, Camilla},
  journal={Mathematics in Engineering},
  volume={5},
  number={1},
  pages={1--41},
  year={2023},
  publisher={Amer Inst Mathematical Sciences-Aims}
}

@article{otto2011rayleigh,
  title={{R}ayleigh-{B}{\'e}nard convection: improved bounds on the Nusselt number},
  author={Otto, Felix and Seis, Christian},
  journal={Journal of Mathematical Physics},
  volume={52},
  number={8},
  pages = {083702},
  year={2011},
  publisher={AIP Publishing}
}

@article{TD2017,
  title={Optimal wall-to-wall transport by incompressible flows},
  author={Tobasco, Ian and Doering, Charles R.},
  journal={Physical Review Letters},
  volume={118},
  number={26},
  pages={264502},
  year={2017},
  publisher={APS}
}

@article{chernyshenko2022relationship,
  title={Relationship between the methods of bounding time averages},
  author={Chernyshenko, Sergei},
  journal={Philosophical Transactions of the Royal Society A},
  volume={380},
  number={2225},
  pages={20210044},
  year={2022},
  publisher={The Royal Society}
}

@article{ChanMal2025,
	author = {Chanillo, Sagun and Malchiodi, Andrea},
	journal = {Inventiones mathematicae},
	number = {2},
	pages = {633--660},
	title = {Sharp bounds on the {N}usselt number in {R}ayleigh-{B}{\'e}nard convection and a bilinear estimate by {C}oifman-{M}eyer},
	volume = {240},
	year = {2025}
}

@book{FL2023,
  title={Stochastic partial differential equations in fluid mechanics},
  author={Flandoli, Franco and Luongo, Eliseo},
  volume={2330},
  year={2023},
  publisher={Springer}
}

@book{CG24,
  title={Stochastic calculus in infinite dimensions and {SPDE}s},
  author={Crisan, Dan and Goodair, Daniel},
  year={2024},
  publisher={Springer}
}

@article{seis2015scaling,
  title={Scaling bounds on dissipation in turbulent flows},
  author={Seis, Christian},
  journal={Journal of Fluid Mechanics},
  volume={777},
  pages={591--603},
  year={2015},
  publisher={Cambridge University Press}
}

@article{lohse2024ultimate,
  title={Ultimate {R}ayleigh-{B}{\'e}nard turbulence},
  author={Lohse, Detlef and Shishkina, Olga},
  journal={Reviews of Modern Physics},
  volume={96},
  number={3},
  pages={035001},
  year={2024},
  publisher={APS}
}

@article{shishkina2024ultimate,
  title={Ultimate regime of {R}ayleigh-{B}{\'e}nard turbulence: {S}ubregimes and their scaling relations for the {N}usselt vs {R}ayleigh and {P}randtl numbers},
  author={Shishkina, Olga and Lohse, Detlef},
  journal={Physical Review Letters},
  volume={133},
  number={14},
  pages={144001},
  year={2024},
  publisher={APS}
}

@article{cifani2026stochastic,
  title={Stochastic modeling of {F}ourier modes in two-dimensional turbulence via filtered white noise},
  author={Cifani, Paolo and Flandoli, Franco and Zanoni, Andrea},
  journal={arXiv preprint arXiv:2605.13671},
  year={2026}
}

@article{kazantsev1968enhancement,
  title={Enhancement of a magnetic field by a conducting fluid},
  author={Kazantsev, Alexander P.},
  journal={Sov. Phys. JETP},
  volume={26},
  number={5},
  pages={1031--1034},
  year={1968}
}

@article{kraichnan1994anomalous,
  title={Anomalous scaling of a randomly advected passive scalar},
  author={Kraichnan, Robert H.},
  journal={Physical Review Letters},
  volume={72},
  number={7},
  pages={1016},
  year={1994},
  publisher={APS}
}

@article{picardo2023polymers,
  title={Polymers in turbulence: stretching statistics and the role of extreme strain rate fluctuations},
  author={Picardo, Jason R. and Plan, Emmanuel L.C. VI M. and Vincenzi, Dario},
  journal={Journal of Fluid Mechanics},
  volume={969},
  pages={A24},
  year={2023},
  publisher={Cambridge University Press}
}

@article{vlad1998diffusion,
  title={Diffusion with intrinsic trapping in two-dimensional incompressible stochastic velocity fields},
  author={Vlad, Mǎdǎlina and Spineanu, Florin and Misguich, Jacques H. and Balescu, Radu},
  journal={Physical Review E},
  volume={58},
  number={6},
  pages={7359},
  year={1998},
  publisher={APS}
}

@article{FLN2026a,
  title={Wall-to-wall heat transport by stochastic velocity fields},
  author={Flandoli, Franco and Lange, Theresa and Nobili, Camilla},
  journal={preprint},
  year={2026}
}

@article{FLN2026c,
  title={Scaling {B}ounds for {S}tochastic {P}orous {M}edia {C}onvection},
  author={Flandoli, Franco and Lange, Theresa and Nobili, Camilla},
  journal={preprint},
  year={2026}
}

@article{doering2019optimal,
  title={On the optimal design of wall-to-wall heat transport},
  author={Doering, Charles R. and Tobasco, Ian},
  journal={Communications on Pure and Applied Mathematics},
  volume={72},
  number={11},
  pages={2385--2448},
  year={2019},
  publisher={Wiley Online Library}
}

@article{howard1963heat,
  title={Heat transport by turbulent convection},
  author={Howard, Louis N.},
  journal={Journal of Fluid Mechanics},
  volume={17},
  number={3},
  pages={405--432},
  year={1963},
  publisher={Cambridge University Press}
}

@article{doering1998bounds,
  title={Bounds for heat transport in a porous layer},
  author={Doering, Charles R. and Constantin, Peter},
  journal={Journal of Fluid Mechanics},
  volume={376},
  pages={263--296},
  year={1998},
  publisher={Cambridge University Press}
}

\end{document}